\documentclass[11pt]{amsart}
\usepackage{amssymb}
\usepackage{amsmath, amsthm, amsopn, amsfonts}
\usepackage[dvips]{graphicx}
\usepackage{graphpap}
\def\Re{\textrm{Re}} 
\def\R{\mathbb R}
\def\Z{\mathbb Z}
\def\T{\mathbb T}
\def\C{\mathbb C}
\def\N{\mathbb N}
\def\11{{\rm 1~\hspace{-1.4ex}l} }

\newtheorem{proposition}{Proposition}
\newtheorem{lemme}[proposition]{Lemma}

\newtheorem{definition}[proposition]{Definition}
\newtheorem{remarque}[proposition]{Remark}
\newtheorem{theoreme}[proposition]{Theorem}

\numberwithin{equation}{section}
\numberwithin{proposition}{section}
\begin{document}
\title[Ill-posedness issues for nonlinear dispersive equations]
{Ill-posedness issues for nonlinear dispersive equations }
\author{N. Tzvetkov}
\address{D\'epartement de Math\'ematiques, Universit\'e Lille I, 59 655 Villeneuve d'Ascq Cedex, France}
\email{nikolay.tzvetkov@math.univ-lille1.fr}
\urladdr{http://math.univ-lille1.fr/~tzvetkov}
\begin{abstract}
These notes are devoted to the notion of well-posedness of the Cauchy problem for nonlinear dispersive equations.
We present recent methods for proving ill-posedness type results for dispersive PDE's.
The common feature in the analysis is that the proof of such results requires the construction of high frequency
approximate solutions on small time intervals (possibly depending on the  frequency).
The classical notion of well-posedness, going back to Hadamard, requires the existence, the uniqueness and the continuity
of the flow map on the spaces where the existence is established. It turns out that in many cases a stronger form
of well-posedness holds. Namely, the flow map enjoys better continuity properties as for example being Lipschitz continuous 
on bounded sets. In such a situation we say that the corresponding problem is semi-linearly well-posed in the corresponding
functional setting. Our main message is that for dispersive PDE's, contrary to the case of hyperbolic PDE's, 
the verification whether an equation in hand is semi-linearly well-posed in a given functional framework requires a considerable care.
Our examples are KdV type equations and non linear Schr\"odinger equations. 
\end{abstract}
\maketitle
\section{Introduction}
We will discuss here the Cauchy problem for nonlinear PDE's which can be written in the form
\begin{eqnarray}\label{1}
\dot{u}(t)=Lu(t)+F(u(t)),\quad u(0)=u_0,
\end{eqnarray}
where $u(t)$, $t\in\R$ is a function defined on a Riemannian manifold $(M,g)$ with values either in $\R$ or in $\C$.
In (\ref{1}), $L$ is a linear map such that $\exp(tL)$ is well-defined and
acts as an isometry on the Sobolev spaces $H^{s}(M)$ while $F(u(t))$ represents the nonlinear interaction. 
The initial data $u_0$ is supposed to belong to $H^{s}(M)$.
This choice is natural because, for the models we are interested in, the equation (\ref{1}) enjoys conservation laws providing a uniform
control on (low regularity) Sobolev norms of the solutions of (\ref{1}).
An important aspect of the analysis of the Cauchy problem (\ref{1}) is to understand the impact of the interplay between $L$ 
and $F$ on the behavior on the solutions of (\ref{1}). Here we will study this issue only for {\it small times} $t$.
As far as the Sobolev spaces $H^{s}(M)$ are chosen for phase spaces, the local in time behavior of the solutions
is naturally linked to the notion of well-posedness of the Cauchy problem
(\ref{1}). Let us now state the notion of well-posedness that will be used here.
\begin{definition}\label{def1}
We say that the Cauchy problem (\ref{1}) is well-posed for data in $H^{s}(M)$, if for every bounded set $B$ of $H^{s}(M)$ 
there exist $T>0$ and a Banach space ${\mathcal X}_{T}$ continuously embedded in $C([-T,T]\,;\,H^{s}(M))$ such that if $u_0\in B$ then 
there exists a unique solution $u$ of (\ref{1}) on $[-T,T]$ in the class ${\mathcal X}_{T}$. 
Moreover :
\begin{enumerate}
\item
The flow map $u_0\mapsto u$ is continuous from $B$ to $C([-T,T]\,;\,H^{s}(M))$.
\item
Higher smoothness is propagated by the flow. More precisely,
if $u_0\in H^{\sigma}(M)$, $\sigma\geq s$ then $u\in C([-T,T]\,;\,H^{\sigma}(M))$.
\end{enumerate}
\end{definition}
Let us notice that in the above definition, the time of existence $T$ depends only on the bounded set $B$, i.e. on an $H^s$
bound of the initial data. There are several important examples of the so called critical problems when the time of existence existence
is depending in a more complicated way on the initial data. It is worth
noticing that ``usually'', if a problem in hand is critical for data
in $H^s$ then it is well-posed in the sense of Definition \ref{def1} for data in $H^{\sigma}$, $\sigma>s$.
It is also ``usual'' that the well-posedness in $H^s$, implies the well-posedness in $H^{s'}$, $s'\geq s$.
In view of the propagation of regularity property~(2) in
Definition~\ref{def1}, on may see the solutions of (\ref{1}) as a limit of
smooth solutions. This may motivate one (see e.g. \cite{KT1,KT2,IK}) to
restrict the study of (\ref{1}) for low regularity data to such solutions
which are limits, in the low regularity topology, of smooth
solutions. Uniqueness in the class of the limits of smooth solutions then
follows from the existence and the propagation of the regularity, and thus
one avoids the difficulty to find a space ${\mathcal X}_{T}$, continuously embedded in $C([-T,T]\,;\,H^{s}(M))$
where the uniqueness holds.
\\

A very common way to prove the well-posedness of (\ref{1}) is to solve by a contraction principle an equivalent integral equation,
exactly as we do in the proof of the Cauchy-Lipschitz theorem in the theory of the ordinary differential equations.
More precisely, the problem (\ref{1}) can be rewritten, at least formally, as an integral equation (Duhamel formula)
\begin{eqnarray}\label{2}
u(t)=\exp(tL)u_{0}+\int_{0}^{t}\exp((t-t')L)F(u(t'))dt' \, .
\end{eqnarray}
The well-posedness of (\ref{1}) is reduced to finding a functional spaces ${\mathcal X}_{\tau}$, $\tau>0$ 
continuously embedded in $C([-\tau,\tau]\,;\,H^{s}(M))$ such that for every  bounded set $B$ of $H^{s}(M)$ there exists $T>0$
such that for every $u_0\in B$ the right hand-side of (\ref{2}) is a contraction in a suitable ball of ${\mathcal X}_{T}$.
In some cases, the space $C([-\tau,\tau]\,;\,H^{s}(M))$ can give the contraction properties. However, in these cases
the assumption on $s$ is quite restrictive. In order to include larger possible values of $s$,
the whole difficulty in making work the above approach is to find functional spaces ${\mathcal X}_{\tau}$, $\tau>0$ which are 
adapted in the best way to the equation in hand. This problematic has now a long history and remains and active research field.
Once the existence and the uniqueness in ${\mathcal X}_{T}$ is established, it is natural to look for a larger uniqueness class,
for instance one may ask whether the uniqueness holds in $C([-T,T]\,;\,H^{s}(M))$ (cf. e.g.~\cite{Kato3}).

It turns out that if we are able to show the well-posedness of (\ref{1}) by the above procedure then the flow map enjoys
better continuity properties, for example it is Lipschitz continuous on $B$, and, in the case of polynomial nonlinearities
it is a $C^{\infty}$ map from $H^{s}(M)$ to $C([-T,T]\,;\,H^{s}(M))$. 
These properties seem to be related to what we call a semi-linearly well-posed problem.
The following definition seems to be natural (cf. e.g. \cite{BKPSV,BPS,Bo5} ...).
\begin{definition}\label{def2}
We say that the Cauchy problem (\ref{1}) is semi-linearly well-posed for data in $H^{s}(M)$, if it is well-posed
in the sense of Definition \ref{def1}, and, in addition the flow map $u_0\mapsto u$ is {\bf uniformly continuous} 
from $B$ to $C([-T,T]\,;\,H^{s}(M))$.
\end{definition}  
The notion of well-posedness of Definition \ref{def1} is invariant under changes of variables in the phase space which are continuous 
on $H^s$. Similarly the notion of semi-linear well-posedness is invariant under uniformly continuous changes of variables. Therefore,
it is not excluded that, by a change of variables (gauge transform)
$$
u(t)\longrightarrow v(t)
$$
in (\ref{1}) which is continuous on $H^s$ but not uniformly continuous, the equation for $v(t)$ to be semi-linearly well-posed even
if the equation for $u(t)$ is not semi-linearly well-posed.
\\

Another and quite different way to solve (\ref{1}) is to apply a compactness argument. Roughly speaking, it means to solve the equation
by passing to a (weak) limit in a family of approximate solutions. Usually this method can provide the well-posedness
of (\ref{1}), but it does not give directly the semi-linear well-posedness as the contraction method does.
\\

A natural question is whether there exists PDE's which are well-posed but not semi-linearly well-posed in $H^s(M)$.
Probably the simplest example of such a PDE is the Burgers equation
\begin{equation}\label{3}
u_t+uu_x=0,
\end{equation}
posed on $H^{s}(\R)$ for real valued $u$ (if $u$ is not real valued the situation is quite different, 
as it is shown in \cite{Christ}).
It turns out that (\ref{3}) 
is well-posed in $H^s(\R)$, $s>3/2$ but not semi-linearly well-posed in this same space.
Let us explain how we prove the well-posedness of (\ref{3}) for data in $H^s(\R)$, $s>3/2$.
Let $u$ be a smooth solution of (\ref{3}) which describes a continuous curve in all $H^{\sigma}$, $\sigma\in\R$.
Our purpose is to establish a priori bounds for $u$.
Denote by $D^s$ the Fourier multiplier with symbol $(1+\xi^2)^{s/2}$, i.e. 
$$
\widehat{D^{s}u}(\xi)=(1+\xi^2)^{s/2}\widehat{u}(\xi)\, ,
$$
where the Fourier transform is defined as follows
$$
\widehat{u}(\xi)=\int_{-\infty}^{\infty}e^{-ix\xi}\, u(x)dx\, .
$$
Notice that $\|u\|_{H^s}=\|D^{s}u\|_{L^2}$.
Applying $D^s$ to (\ref{3}), multiplying it with $D^s u$ and an integration by parts gives,
$$
\frac{d}{dt}\|u(t,\cdot)\|_{H^s}^{2}
=
\int_{-\infty}^{\infty}u_{x}(t,x)\big(D^{s}u(t,x)\big)^{2}dx
-2\int_{-\infty}^{\infty}\big([D^s,u]u_{x}\big)(t,x)D^{s}u(t,x)dx \, .
$$
Using the Kato-Ponce (cf. \cite{KP}) commutator estimate
\begin{equation}\label{Kato-Ponce}
\|[D^s,f]\, g\|_{L^2}\leq C\big(\|f_x\|_{L^{\infty}}\|D^{s-1}g\|_{L^2}+\|D^s f\|_{L^2}\|g\|_{L^{\infty}}\big)
\end{equation}
with $f=u$ and $g=u_x$, we obtain that
$$
\frac{d}{dt}\|u(t,\cdot)\|_{H^s}^{2}\leq C\|u_{x}(t,\cdot)\|_{L^{\infty}}\|u(t,\cdot)\|_{H^s}^{2}\,\, .
$$
Thus the Gronwall lemma yields that for every $0\leq t\leq T$,
\begin{equation}\label{KAM1}
\|u(t,\cdot)\|_{H^s}\leq \|u(0,\cdot)\|_{H^s}\,\exp\big(C\|u_x\|_{L^{1}([0,T]\,;\,L^{\infty})}\big)\,. 
\end{equation}
If $s>3/2$, the Sobolev embedding gives,
\begin{equation}\label{KAM2}
\|u_x\|_{L^{1}([0,T]\,;\,L^{\infty})}\leq C\,T\|u\|_{L^{\infty}([0,T]\,;\,H^{s})}\, .
\end{equation}
Combining (\ref{KAM1}) and (\ref{KAM2}), using a continuity argument, we deduce that there exist $c>0$ and $C>0$ such that if
$$
T\leq c(1+\|u_0\|_{H^s})^{-1}
$$
then
\begin{equation}\label{key}
\|u_x\|_{L^{1}([0,T]\,;\,L^{\infty})}\leq C
\end{equation}
and 
\begin{equation}\label{apriori}
\|u\|_{L^{\infty}([0,T]\,;\,H^{s})}\leq C\|u(0,\cdot)\|_{H^s} \,.
\end{equation}
The priori estimate (\ref{apriori}) is the key to perform a classical compactness argument (cf. e.g. \cite{Lions}) 
which provides the existence. More precisely one passes into the limit
$\varepsilon\rightarrow 0^{+}$ in the solutions of a regularized equation, e.g.
$$
u^{\varepsilon}_{t}-\varepsilon
u^{\varepsilon}_{txx}+u^{\varepsilon}u^{\varepsilon}_{x}=0 \quad {\rm (BBM\,\, regularization) }
$$
or
$$
u^{\varepsilon}_{t}-\varepsilon
u^{\varepsilon}_{xxxx}+u^{\varepsilon}u^{\varepsilon}_{x}=0 \quad {\rm (parabolic\,\, regularization) }
$$
etc. The particular choice of the regularized equation is not of importance,
the crucial point in that exactly as above one may show that the solutions of
the regularized equation enjoy the bounds (\ref{key}) and (\ref{apriori})
{\it uniformly in} $\varepsilon$ which enables one to passe into the limit
$\varepsilon\rightarrow 0^{+}$.
\\

The uniqueness is easily ensured by the Gronwall lemma and and the control on
$u_x$ in $L^{\infty}$.
The propagation of the higher Sobolev regularity readily follows from (\ref{key}) and (\ref{apriori}).
\\

The continuous dependence is a slightly more delicate issue and can be obtained for instance by the Bona-Smith argument \cite{BS}
(cf. also \cite{Kato1}). 
Let us briefly recall this argument. Fix a bump function $\rho\in{\mathcal S}(\R)$ such that 
$\widehat{\rho}\in C_{0}^{\infty}(\R)$ and $\widehat{\rho}(\xi)=1$ for $\xi$ 
in a neighborhood of $0$. For $\varepsilon>0$, we set
$
\rho_{\varepsilon}(x):=\varepsilon^{-1}\rho(x/\varepsilon)\, .
$
Let $u$ be a solution of (\ref{3}) with data $u(0)\in H^{s}(\R)$, $s>3/2$ which belongs to {\it a fixed bounded set} of $H^s(\R)$.
Denote by $u^{\varepsilon}$ the solution of the Burgers equation (\ref{3}) with initial data $\rho_{\varepsilon}\star u(0)$.
One can easily check that
$$
\|\rho_{\varepsilon}\star u(0)\|_{H^s}\leq C\|u_0\|_{H^s},\quad \varepsilon\in ]0,1]
$$
and therefore we can assume that $u^{\varepsilon}$ enjoys the bounds (\ref{key}) and (\ref{apriori}) on the time of existence of $u$.
For $\varepsilon>\varepsilon'>0$, we set $v:=u^{\varepsilon}-u^{\varepsilon'}$. Then $v$ is a solution of the equation
\begin{equation}\label{sol1}
v_t-vv_x+u^{\varepsilon}_{x}v+u^{\varepsilon}v_x=0\, .
\end{equation}
Note that we privilege $u^{\varepsilon}$ to $u^{\varepsilon'}$ because $\varepsilon>\varepsilon'$.
It is easy to check that
\begin{equation}\label{sol2}
\|v(0)\|_{H^s}=o(1),\quad \|v(0)\|_{L^2}=o(\varepsilon^{s})
\end{equation}
as $\varepsilon\rightarrow 0$. 
Multiplying (\ref{sol1}) with $v$ and applying (\ref{key}) (with
$u^{\varepsilon}$ and $u^{\varepsilon'}$ instead of $u$) gives the bound
\begin{equation}\label{sol2bis}
\|v(t,\cdot)\|_{L^2}\leq o(\varepsilon^{s})
\end{equation}
for $t$ in the time of existence of $u$.
Applying $D^s$ to (\ref{sol1}), multiplying it with $D^s v$ and using the Kato-Ponce estimate (\ref{Kato-Ponce}) yields the estimate
\begin{multline}\label{kak}
\frac{d}{dt}\|v(t,\cdot)\|_{H^s}^{2}
\leq
C\|v(t,\cdot)\|_{H^s}^{3}
+
C\|u^{\varepsilon}(t,\cdot)\|_{H^{s}}\|v(t,\cdot)\|_{H^s}^{2}
+
\\
+
C\|u^{\varepsilon}(t,\cdot)\|_{H^{s+1}}\|v(t,\cdot)\|_{H^{s-1}}\, \|v(t,\cdot)\|_{H^s}.
\end{multline}
The third term in the right hand-side of (\ref{kak}) is a new one compared to
the a priori bound discussion above. We may have that 
$\|u(t,\cdot)\|_{H^{s+1}}$ equals infinity but for $\varepsilon > 0$ the quantity
$\|u^{\varepsilon}(t,\cdot)\|_{H^{s+1}}$ exists but is probably very big.
More precisley, using (\ref{apriori}) (with $s+1$ instead of $s$) gives
\begin{equation}\label{patrick1}
\|u^{\varepsilon}(t,\cdot)\|_{H^{s+1}}\leq C \|u^{\varepsilon}(0,\cdot)\|_{H^{s+1}}\leq C\varepsilon^{-1}
\, .
\end{equation}
On the other hand, thanks to (\ref{sol2bis}),
\begin{equation}\label{patrick2}
\|v(t,\cdot)\|_{H^{s-1}}\leq \|v(t,\cdot)\|_{L^2}^{\frac{1}{s}}\|v(t,\cdot)\|_{H^s}^{1-\frac{1}{s}}
\leq o(\varepsilon) \|v(t,\cdot)\|_{H^s}^{1-\frac{1}{s}}
\end{equation}
Using (\ref{kak}), (\ref{patrick1}), (\ref{patrick2}), the Gronwall lemma, and (\ref{sol2}) gives
\begin{equation}\label{sol3}
\|u^{\varepsilon}(t,\cdot)-u^{\varepsilon'}(t,\cdot)\|_{H^s}=\|v(t,\cdot)\|_{H^s}
\leq
C\|u^{\varepsilon}(0,\cdot)-u^{\varepsilon'}(0,\cdot)\|_{H^s}
+o(1)
=
o(1)
\end{equation}
as $\varepsilon\rightarrow 0$.
We can now easily obtain the continuity of the flow map. Indeed, let $(u_{0,n})$ be a sequence converging to $u_0$
in $H^s(\R)$ with corresponding solutions $(u_n)$. 
Using the triangle inequality, we can write
\begin{multline*}
\|u(t,\cdot)-u_n(t,\cdot)\|_{H^s}\leq 
\|u(t,\cdot)-u^{\varepsilon}(t,\cdot)\|_{H^s}
\\
+
\|u^{\varepsilon}(t,\cdot)-u^{\varepsilon}_{n}(t,\cdot)\|_{H^s}
\\
+
\|u^{\varepsilon}_n(t,\cdot)-u_{n}(t,\cdot)\|_{H^s}
:=(1)+(2)+(3)\, .
\end{multline*}
Using that
$ 
\rho_{\varepsilon}\star u_{0,n}\rightarrow u_{0,n}
$
as $\varepsilon\rightarrow 0$
in $H^s(\R)$, {\it uniformly in $n$}, 
exactly\footnote{The important point is that each of the terms $(1)$, $(2)$, $(3)$ contains
a solution involving $\varepsilon$. }
as in the proof of (\ref{sol3}), we can show that as $\varepsilon\rightarrow 0$,
$$
(1)\leq o(1),\quad (2)\leq C\|u_{0,n}-u_0\|_{H^s}+o(1),\quad (3)\leq o(1)
$$
which clearly implies the continuity on $H^s(\R)$, $s>3/2$ of the flow map of the Burgers equation (\ref{3}).
\\

At this point, it is worth to notice that the argument based on a priori estimates for proving the well-posedness that 
we have just presented is less perturbative (``more nonlinear'') than the contraction method explained after Definition~\ref{def1}.
It has the advantage to have a larger scope of applicability compared to the contraction method, but, at the present moment,
to make it work one should require considerably more regularity on the initial data. 
\\
 
Let us next describe an argument providing the lack of semi-linear well-posedness of~(\ref{3}).
We first observe that if $u$ solves (\ref{3}) then so does 
\begin{equation}\label{simetria}
v(t,x)= u(t,x-\omega t) + \omega,\quad \omega\in\R.
\end{equation}
The shift in the spatial variable in (\ref{simetria}) is ``responsible'' for the failure of uniform continuity of the flow map.
The constant $\omega$ in (\ref{simetria}) can be replaced by a function which is zero at infinity, thanks to the
finite propagation speed of the Burgers equation. More precisely, inspired by (\ref{simetria}), we look for an approximate 
solution of the Burgers equation of the form
\begin{equation}\label{form}
u_{ap}^{\omega,\lambda}(t,x)=\omega\lambda^{-1}\widetilde{\varphi}\big(x/\lambda^{\delta}\big)+
\lambda^{-\delta/2-s}\varphi\big(x/\lambda^{\delta}\big)\cos(\lambda x-\omega t),
\end{equation}
where $s>3/2$, $t\in [-1,1]$, $\omega\in\R$, $1<\delta<2$ and $\varphi$, $\widetilde{\varphi}$ are non zero $C_{0}^{\infty}(\R)$
functions such that $\widetilde{\varphi}$ is equal to one on the support  $\varphi$. 
We can then show that there exists $\varepsilon>0$ such that
\begin{equation}\label{form-bis}
\|\partial_{t}u_{ap}^{\omega,\lambda}+u_{ap}^{\omega,\lambda}\partial_{x}u_{ap}^{\omega,\lambda}\|_{L^2(\R)}
\leq
C\, \lambda^{-\varepsilon-s}\, .
\end{equation}
Thanks to  (\ref{form-bis}) and the well-posedness analysis in $H^s$, $s>3/2$, we obtain that (\ref{form}) is indeed
a good approximate solution, in $H^s$, of the Burgers equation. Considering the sequences $(u_{ap}^{1,\lambda})$
and $(u_{ap}^{-1,\lambda})$, $\lambda\gg 1$ gives the failure of uniform continuity on $H^s(\R)$, $s>3/2$ 
of the flow map of the Burgers equation. 
\\

However, the Burgers equation (\ref{3}) does not fit in the class of dispersive PDE's and one may think that the above 
described property of (\ref{3}) is only related to its hyperbolic nature.
It turns out that the Benjamin-Ono equation
\begin{equation}\label{4}
u_t+Hu_{xx}+uu_x=0,
\end{equation}
posed on $H^{s}(\R)$ (in (\ref{4}) $H$ denotes the Hilbert transform which is a ``zero order'' operator) 
is well-posed in $H^s(\R)$, $s>3/2$ but not semi-linearly well-posed in this same space.
The equation (\ref{4}) fits in the class of the dispersive equations because of the presence of the term $Hu_{xx}$.
We use the term ``dispersive equation'' since any solution of
\begin{equation}\label{5}
u_t+Hu_{xx}=0,
\end{equation}
issued form $L^{1}(\R)$ initial data disperses as $t\rightarrow\infty$, more precisely, 
$$
\lim_{t\rightarrow\infty}\|u(t,\cdot)\|_{L^{\infty}(\R)}=0\, .
$$
However this property is for large times, and, since we are concerned with a small time analysis a more relevant property related to
the dispersive nature of the 
equation (\ref{5}) is the (small time) Strichartz inequality (cf. e.g. \cite{Ponce}). 
More precisely, there exists $C>0$ such that for every $T>0$, 
every $u_0\in L^{2}(\R)$
the solution $u$ of (\ref{5}) with data $u_0$ satisfies,
\begin{equation}\label{6}
\|u\|_{L^{p}([0,T]\,;\,L^{q}(\R))}\leq C\|u_0\|_{L^2(\R)},\quad \frac{2}{p}+\frac{1}{q}=\frac{1}{2},\,\,\, p\geq 4\, .
\end{equation}
Estimates of type (\ref{6}) are usually very useful to apply the contraction strategy but in the case of (\ref{4})
they are not sufficient to make it work.
\\

We next consider the KdV equation 
\begin{equation}\label{7}
u_t+u_{xxx}+uu_x=0,
\end{equation}
posed on $H^{s}(\R)$, which have a higher order dispersion compared to (\ref{4}).
It turns out that, in sharp contrast with the Burgers and the Benjamin-Ono equations, 
the KdV equation  (\ref{7}) is semi-linearly well-posed for data in $H^{s}(\R)$, $s>-3/4$.
Therefore, the notion of semi-linear well-posedness makes a natural classification in the class of KdV type models, i.e.
equation (\ref{1}) with $F(u)=uu_x$ and $L=|D_{x}|^{\alpha}\partial_x$, depending on the order of dispersion $\alpha>0$.
\\

Another set of models where the notion of semi-linear well-posedness is naturally involved (but less understood) are
the nonlinear Schr\"odinger equations (NLS). Let $(M,g)$ be a compact smooth boundaryless Riemannian manifold of dimensions~$d=2,3$. 
Denote by $\Delta$ the Laplace-Beltrami operator associated to the metric $g$.
We consider the nonlinear Schr\"odinger equation
\begin{equation}\label{8}
iu_{t}+\Delta u-|u|^{2}u=0,
\end{equation}
posed on $M$. In (\ref{8}), $u$ is complex valued function on $M$. 
Let us first consider the case $d=2$, i.e. the case when $M$ is a surface.
If $M$ is the flat torus $\T^2$ then the Cauchy problem associated to (\ref{8}) is semi-linearly well-posed for data
in $H^{s}(\T^2)$, provided $s>0$. This result is essentially sharp, since for $s<0$ the problem is not semi-linearly
well-posed for data in in $H^{s}(M)$ for an arbitrary $(M,g)$. On the other hand if $M$ is the standard sphere $S^2$ then 
the Cauchy problem for (\ref{8}) is not semi-linearly well-posed far data in $H^s(S^2)$, $s<1/4$, in sharp contrast
with the case of the torus $\T^2$. Hence the same equation (\ref{8}) behaves quite differently with respect to the
semi-linear well-posedness depending on the geometry of the spatial domain.
It is an interesting open problem whether (\ref{8}) posed on $S^2$ might be well-posed for some $s\in[0,1/4]$.
Let us next consider the case $d=3$. If $M$ is the torus $\T^3$ or the sphere $S^3$ then the Cauchy problem for (\ref{8})
is semi-linearly well-posed for data in $H^s$, $s>1/2$. It turns out that this result is essentially 
sharp even regarding the classical notion of well-posedness. More precisely, for $0<s<1/2$, the Cauchy problem for 
(\ref{8}) posed on an arbitrary $M$ is not well-posed for data in $H^s$.
\\

One may ask for the critical threshold in the scale of $H^s$ for the well-posedness  of (\ref{1}). 
It means to find a real number $s_c$ such that for $s<s_c$ (\ref{1}) is not well-posed for data in $H^s(M)$, 
while for $s>s_c$ (\ref{1}) is well-posed for data in $H^s(M)$.
Similarly, one can define a critical threshold for the semi-linear well-posedness.
In this context, the discussion around (\ref{8}) above simply affirms that, for $d=3$ and $M=\T^3$, the value $s_c=\frac{1}{2}$
is the critical threshold for both the well-posedness and the semi-linear well-posedness, as far as positive values of the
Sobolev regularity $s$ are considered.
It is a natural question whether the  critical threshold for the well-posedness and the semi-linear well-posedness may be different. 
The answer of this question is positive as shows the following example. Consider the following version of the modified KdV equation
\begin{equation}\label{8.5}
u_t+u_{xxx}+(u^{2}-\int_{\T}u^{2}(t,y)dy) u_x=0,
\end{equation}
posed on the torus $\T=\R/ \Z$.
The equation (\ref{8.5}) can be obtained from the modified KdV equation
\begin{equation}\label{9}
v_t+v_{xxx}+v^{2}v_x=0,
\end{equation}
by the {\it gauge transformation} $u\rightarrow v$ defined as
$$
v(t,x)=u\big(t,x-\int_{0}^{t}\int_{\T}u^{2}(\tau,y)dyd\tau\big)\, .
$$
The Cauchy problem for (\ref{8.5}) is semi-linearly well-posed for data in $H^{s}(\T)$, $s>1/2$ (cf. \cite{Bo3}), it is not
semi-linearly well-posed for data in $H^{s}(\T)$, $3/8<s<1/2$ (cf. \cite{TaTs}), but ... it is still well-posed for data in 
$H^{s}(\T)$, $s\in [3/8,1/2]$ (cf. \cite{TaTs,KT1,KT2}).
Hence the  critical threshold for the well-posedness and the semi-linear well-posedness can be different.
\\

Let us complete this introduction by noticing that, in the last years, gauge transformations were an important tool in
the study of dispersive PDE's,
cf. e.g. \cite{Bo3,HO,MR,TaTs,Tao1,Tao2,Tao3}. In the light of the above discussion,
one may wish to see the gauge  transformations as a tool which, essentially
speaking, transforms a problem which is not  semi-linearly well-posed to a
problem which is semi-linearly well-posed. 
\\

{\bf Acknowledgment.} These notes are based on a course given in University
of Hokkaido in September 2004 under an invitation of Professor T. Ozawa and 
Professor Y. Tsutsumi. It is a pleasure to thank all participants for
their interest in the subject of the lectures.   
\section{KdV type problems}
Consider the Cauchy problem for the Korteweg de Vries (KdV) equation
\begin{equation}\label{2.1}
u_{t}+u_{xxx}+uu_{x}=0,\quad u(0)=u_0.
\end{equation}
The best known result regarding the well-posedness of (\ref{2.1}) is due to Kenig-Ponce-Vega.
\begin{theoreme}[cf. \cite{KPV3}]\label{thm1}
For $s>-3/4$ the Cauchy problem (\ref{2.1}) is semi-linearly well-posed for data in $H^{s}(\R)$.
\end{theoreme}
To prove Theorem \ref{thm1} one uses the contraction method as explained after Definition \ref{def1} of the previous section.
The spaces $X_T$ where one performs the argument are the Fourier transform restriction spaces introduced by Bourgain \cite{Bo1,Bo2,Bo3}, 
equipped with the norm
$$
\|u\|_{X_{T}}=\inf\{\|w\|_{X},\quad  w\in X \quad{\rm with }\quad w|_{[-T,T]}=u\},
$$
where
$$
\|w\|_{X}^{2}=\int_{\R^2}
\big(1+|\tau-\xi^{3}|^{2}\big)^{b}\big(1+|\xi|^{2}\big)^{s}|\widehat{w}(\tau,\xi)|^{2}d\tau d\xi
$$
with $b>1/2$ sufficiently close to $1/2$. 
The spaces of Bourgain are very useful to recover the derivative loss in the nonlinearity.
We refer to \cite{Bo-book,Gi} for an introduction to the Fourier transform restriction method of Bourgain.
There has been a number of works preceding Theorem \ref{thm1} where the well-posedness
for bigger values of $s$ were established (cf. e.g. \cite{Saut,KPV1,KPV2,Bo3}).
A particularly important step was done in \cite{KPV2}, where it is realized for the first time that the KdV equation can be semi-linearly
well-posed. The value $s=-3/4$ in Theorem~\ref{thm1} is optimal, as far as the semi-linear well-posedness is concerned (cf. \cite{CCT1}).
But it is a priori not excluded (\ref{2.1}) to be well-posed for some $s<-3/4$. 

Next, we consider the Cauchy problem for the Benjamin-Ono (BO) equation (cf.~\cite{Be})
\begin{equation}\label{2.2}
u_{t}+Hu_{xx}+uu_{x}=0,\quad u(0)=u_0.
\end{equation} 
In (\ref{2.2}), $H$ denotes the Hilbert transform, namely,
$$
(Hf)(x):=2\lim_{\varepsilon\rightarrow 0}\int_{|x-y|\geq \varepsilon}\frac{f(y)}{x-y}dy\, .
$$
It is easy to check that for $f\in L^2$,
$$
\widehat{Hf}(\xi)=-i\,{\rm sign}(\xi)\widehat{f}(\xi)\, .
$$
Therefore the Hilbert transform is acting essentially as a zero order operator. The presence of $H$ in (\ref{2.2})
is important to establish some monotonicity properties of the local mass of the solutions of (\ref{2.2}), but it will
not play an essential role in our discussion here.

There has been many works regarding the well-posedness of (\ref{2.2}) (cf. \cite{Saut,ABFS,Iorio,Ponce,KT1,KK,Tao3,BP,IK}).
Let us state a result which is due to Tao.
\begin{theoreme}[cf. \cite{Tao3}]
\label{thm2}
For $s\geq 1$ the Cauchy problem (\ref{2.2}) is well-posed for data in $H^{s}(\R)$.
\end{theoreme}
One may ask whether, similarly to the KdV case, we also have the semi-linear well-posedness in Theorem \ref{thm2}.
It turns out that the answer is negative.
\begin{theoreme}
[cf. \cite{KTz2}]
\label{thm3}
In Theorem \ref{thm2}, one can not replace the well-posedness with semi-linear well-posedness.
\end{theoreme}
Therefore in the well-posedness analysis of (\ref{2.2}), it is not a question to find a suitable space to perform the
contraction method, simply this method for proving the well-posedness does not work, as far as the classical Sobolev spaces
$H^s$ are considered as a space for the initial data.
This fact was first detected in \cite{MST1}.

A related to Theorem \ref{thm3} result is obtained in \cite{BL} where it is shown the lack of semi-linear well-posedness for 
(\ref{2.2}) with data in $H^s(\R)$, $s<-1/2$.
\\

Interestingly, the modified Benjamin-Ono equation
\begin{equation*}
u_{t}+Hu_{xx}+u^{2}u_{x}=0,\quad u(0)=u_0
\end{equation*}
turns out to be semi-linearly well-posed for data in $H^s(\R)$, $s>1/2$ (cf. \cite{MR}).
Hence, even if the dispersion is the same, the semi-linear well-posedness may also be sensitive 
to the ``degree'' of the nonlinearity\footnote{The example of KdV and Benjamin-Ono equations is an instance
when the the  semi-linear well-posedness depends on the degree of the dispersion.}.
\\

It is clear that Theorem \ref{thm3} is a consequence of the following statement.
\begin{theoreme}
[cf. \cite{KTz2}]
\label{thm4}
Let $s>0$. There exist two positive constants $c$ and $C$ and 
two sequences $(u_n)$ and $(\widetilde{u}_{n})$ of solutions of the Benjamin-Ono equation 
such that for every $t\in [0,1]$,
$$
\sup_{n}\|u_{n}(t,\cdot)\|_{H^{s}(\R)}+\sup_{n}\|\widetilde{u}_{n}(t,\cdot)\|_{H^{s}(\R)}\leq C\, ,
$$
$(u_n)$ and $(\widetilde{u}_{n})$ satisfy initially
$$
\lim_{n\rightarrow\infty}\|u_{n}(0,\cdot)-\widetilde{u}_{n}(0,\cdot)\|_{H^{s}(\R)}=0,
$$
but, for every $t\in [0,1]$,
$$
\liminf_{n\rightarrow \infty}\|u_{n}(t,\cdot)-\widetilde{u}_{n}(t,\cdot)\|_{H^{s}(\R)}\geq c\,\sin t\, .
$$
\end{theoreme} 
In  the proof of Theorem \ref{thm4}, we will make use of the following well-posedness result for (\ref{2.2}).
\begin{proposition}\label{wp}
Fix $s\geq\sigma>3/2$. Then for every $u_{0}\in H^{s}(\R)$ there exists a unique global
solution $u\in C(\R;H^{s}(\R))$ of (\ref{2.2}). Moreover  
$$
\|u(t,\cdot)\|_{H^{s}(\R)}\leq C\|u_0\|_{H^{s}(\R)}\,\, ,
$$
provided $|t|\leq \min\big(1,c\|u_0\|_{H^{\sigma}(\R)}^{-4}\big)$.
\end{proposition} 
\begin{proof}[Sketch of the proof]
The proof of Proposition~\ref{wp} is based on the compactness argument explained in the introduction in the context of the
Burgers equation (\ref{3}). One first proves the result for $s=\sigma$. 
The nature of the restriction $|t|\leq  c\|u_0\|_{H^{\sigma}(\R)}^{-4}$ is related to the scaling of (\ref{2.2}). 
It turns out that one can reduce the matters to the problem of existence
on the time interval $[0,1]$ with initial data with 
small norm in $H^{\sigma}(\R)$.
Suppose that there  exists a positive constant $\gamma$ such that if the initial data of the Benjamin-Ono equation satisfies
$\|u_0\|_{H^{\sigma}}\leq \gamma$
then we can find a unique solution on the time interval $[0,1]$. We now prove that for $u_0\in H^{\sigma}$ of arbitrary size we can solve
(\ref{2.2}) for time of order $(1+\|u_0\|_{H^{\sigma}})^{-4}$.  Indeed, given $u_{0}\in H^{\sigma}$ we choose $\lambda\ll 1$ such that
\begin{eqnarray}\label{small} 
0<\lambda^{\frac{1}{2}}(1+\lambda^{\sigma})\|u_0\|_{H^{\sigma}}\leq \gamma.
\end{eqnarray}
Set $\widetilde{u}_{0}(x):=\lambda u_{0}(\lambda x)$.  Then due to (\ref{small}),
$
\|\widetilde{u}_{0}\|_{H^{\sigma}}\leq \gamma
$
and we can apply our assumption to $\widetilde{u}_{0}$. Let $\widetilde{u}(t,x)$ be the solution of the Benjamin-Ono
equation with data $\widetilde{u}_{0}$ up to time one. Then one can easily verify that
$
u(t,x):=\lambda^{-1}\widetilde{u}(\lambda^{-2}t,\lambda^{-1}x)
$
is a solution of the Benjamin-Ono equation up to time $\lambda^{2}$ which in view of (\ref{small}) 
is of order $(1+\|u_0\|_{H^{\sigma}})^{-4}$.
Hence we may reduce  the matters to the existence on $[0,1]$ for small data.

Let $u$ be a sufficiently smooth in the scale of the Sobolev spaces solution of the Benjamin-Ono equation (\ref{2.2}).
Then, as in the case of the Burgers equation, one gets the bound
\begin{eqnarray}\label{KP}
\|D^{\sigma}u\|_{L^{\infty}([0,T]\,;\,L^2)}\leq
\|u(0)\|_{H^{\sigma}}\exp(C\|u_x\|_{L^{1}([0,T]\,;\,L^{\infty})})\, .
\end{eqnarray}
Notice that the key quantity $\int_{0}^{T}\|u_{x}(t)\|_{L^{\infty}}dt$ is invariant with respect the scaling of the equation.
More precisely if 
$
u(t,x)=\lambda^{-1}\widetilde{u}(\lambda^{-2}t,\lambda^{-1}x)
$
then
$$
\int_{0}^{\lambda^2}\|u_{x}(t)\|_{L^{\infty}}dt=\int_{0}^{1}\|\widetilde{u}_{x}(t)\|_{L^{\infty}}dt\, .
$$
With
$$
F(T):=\|u_x\|_{L^1_T L^{\infty}}+\|D^{\sigma}u\|_{L^{\infty}_{T}L^{2}},
\quad T\in [0,1]
$$
we can deduce from (\ref{KP}) and the Sobolev inequality (here we use that $\sigma>3/2$) that
$$
F(T)\leq C \|u(0)\|_{H^{\sigma}}\exp(cF(T)).
$$
Now a straightforward continuity argument shows that there exist positive constants $\gamma$ and $C$ such that 
if $\|u(0)\|_{H^{\sigma}}\leq \gamma$ (and hence $F(0)\leq \gamma$) then $F(1)\leq C$, and in particular 
\begin{eqnarray}\label{comp1}
\int_{0}^{1}\|u_{x}(t)\|_{L^{\infty}}dt\leq C
\end{eqnarray}
Using (\ref{comp1}) and (\ref{KP}) (with $T=1$) we obtain that if $u$ is a
smooth solution of the Benjamin-Ono equation, then
\begin{eqnarray}\label{comp2}
\|D^{{\sigma}} u\|_{L^{\infty}([0,1];L^2(\R))}\leq C \|u(0)\|_{H^{\sigma}},
\end{eqnarray}
provided $\|u(0)\|_{H^{\sigma}}\leq \gamma$.
Moreover the solution satisfies (\ref{comp1}).
The bounds (\ref{comp1}) and (\ref{comp2}) enable one to perform a compactness argument for the proof of the existence.
As for the Burgers equation,
the uniqueness follows from the Gronwall lemma, the assumption $\sigma>3/2$ and the Sobolev embedding.
Let us next show the bound for the higher Sobolev norms. Let $s>\sigma$. Then we clearly have an analog of (\ref{KP}) on the 
$H^s$ level. Namely,
\begin{eqnarray*}
\|D^{s}u\|_{L^{\infty}([0,1];L^2)}\leq \|u(0)\|_{H^{s}}\exp(c\int_{0}^{1}\|u_{x}(t)\|_{L^{\infty}}dt)
\leq C\|u(0)\|_{H^{s}},
\end{eqnarray*}
where in the last inequality, we used (\ref{comp1}).
Finally the global well-posedness follows from the conservation lows enjoyed by the solutions of the (\ref{2.2}).
Indeed one has controls (cf. e.g. \cite{ABFS}) on $\|u(t,\cdot)\|_{H^{k/2}}$ for $k=0,1,2,\dots$ Hence the assertion of global existence is
straightforward for $s\geq 2$. For $s<2$ one may use the $H^{3/2}$ well-posedness result result of Theorem~\ref{thm2} and 
the $H^{3/2}$ control. 
This completes the discussion on the proof of Proposition~\ref{wp}.
\end{proof}
Next, we pick a usual bump function $\varphi\in C_{0}^{\infty}(\R)$ such that 
$\varphi(x)=1$ for $|x|<1$ and $\varphi(x)=0$ for $|x|>2$.
Let $\tilde{\varphi}\in C_{0}^{\infty}(\R)$ be equal to one on the support of $\varphi$. Notice that $\varphi\tilde{\varphi}=\varphi$. 
For $0<\delta<1$, we set
\begin{equation}\label{lambda}
\varphi_{\lambda}(x):=\varphi\big(\frac{x}{\lambda^{1+\delta}}\big),\quad
\tilde{\varphi}_{\lambda}(x):=\tilde{\varphi}\big(\frac{x}{\lambda^{1+\delta}}\big)\, .
\end{equation}
The assertion of Theorem~\ref{thm4} is a corollary of the following statement.
\begin{theoreme}\label{thm5}
Let  $\max(1-s,0)<\delta< 1$ and $|\omega| \ll \lambda^{\frac{1-\delta}2}$. Let $u_{\omega,\lambda}$ 
be the unique global solution of the Benjamin-Ono equation subject to initial data
$$
u_{\omega,\lambda}(0,x)
=
-\omega\,\lambda^{-1}\tilde{\varphi}_\lambda(x)
-
\lambda^{-\frac{1}{2}-\frac{\delta}{2}-s} \varphi_\lambda(x)
\cos \lambda x.
$$
Then the identity
$$
u_{\omega,\lambda}(t,x) =-\lambda^{-\frac{1}{2}-\frac{\delta}{2}-s} 
\varphi_\lambda(x)
\cos (-\lambda^2 t + \lambda x + \omega t)  + 
{\mathcal O}\left(\lambda^{- \frac{\min\{\delta, 1-\delta\}}{4(s+2)} } 
+ |\omega|\, \lambda^{-\frac{1-\delta}2}\right)
$$
holds in $H^{s}_{x}(\R)$, uniformly in $t\in[0,1]$ and $\lambda\gg 1$.
\end{theoreme}
Let us notice that if $\omega=0$, the solution propagates as a high frequency linear Benjamin-Ono wave while when $\omega\neq  0$, 
the solution propagates as a high frequency linear dispersive wave with modified propagation speed which is the crucial nonlinear effect.
\\

Let us now show why Theorem~\ref{thm5} implies Theorem~\ref{thm4}. Apply Theorem \ref{thm5} 
with $\omega=\pm 1$ and $\lambda=1,2,\dots$ We thus obtain two families $(u_{1,\lambda})$ and $(u_{-1,\lambda})$
of solutions to  the Benjamin-Ono equation. 
Notice that
$$
\|u_{1,\lambda}(0,\cdot)-u_{-1,\lambda}(0,\cdot)\|_{H^{s}}
\leq C \lambda^{-\frac{1-\delta}{2}} 
$$
and moreover due to Theorem~\ref{thm5}, setting  $\kappa= -\lambda^{2} t+ \lambda x$, we arrive at
$$
\|u_{1,\lambda}(t,\cdot)-u_{-1,\lambda}(t,\cdot)\|_{H^{s}}=
\Vert 
\lambda^{-(\frac{1+\delta}2+s)} \varphi_{\lambda}(x) (\cos(\kappa+t) - \cos(\kappa-t )) \Vert_{H^s_x} + o(1), 
$$
if $t\in[0,1]$ and where $o(1)\rightarrow 0$ as $\lambda\rightarrow\infty$.
At this point we need the following elementary lemma whose proof will be omitted.
\begin{lemme}\label{loc}
Fix $s\geq 0$, $0<\delta<1$, $\alpha\in\R$ and  $\varphi\in C_{0}^{\infty}(\R)$.
Then 
$$
\lim_{\lambda \to \infty} \lambda^{-\frac{1+\delta}{2}-s} 
\left\|\varphi_\lambda(x)  \cos(\lambda x+\alpha)\right\|_{H^{s}_{x}}
= \frac1{\sqrt{2}} \Vert \varphi \Vert_{L^2}\, ,
$$
where $\varphi_{\lambda}$ is defined by (\ref{lambda}).
\end{lemme}
Using Lemma~\ref{loc}, we get
\begin{eqnarray*}
\lim_{\lambda \to \infty}
\Vert \lambda^{-(\frac{1+\delta}2+s)} \varphi_{\lambda}(x) (\cos(\kappa   +t ) - \cos(\kappa-t )) \Vert_{H^s_x}=\sqrt{2}|\sin t|  \,\,
\Vert \varphi \Vert_{L^2}\, .
\end{eqnarray*} 
Therefore Theorem~\ref{thm5} implies Theorem~\ref{thm4}.

\begin{proof}[Proof of Theorem~\ref{thm5}]
Let $u_{low}(t,x)$ be the solution of (\ref{2.2}) with initial data
\begin{equation}\label{low}
u_{low}(0,x)=-\omega\lambda^{-1}\tilde{\varphi}_\lambda(x),\quad 0<\delta<1,
\quad \omega\in\R.
\end{equation}
In the next lemma, we collect several bounds for $u_{low}(t,x)$.
\begin{lemme}\label{ulow}  
Let $k\geq 0$. Then the following estimates hold :
\begin{eqnarray}
\|\partial_x^k u_{low}(t,\cdot)\|_{L^2(\R)} & \leq & C |\omega|\,\lambda^{-\frac{1-\delta}{2}-k(1+\delta)},
\label{parvo}
\\
\|\partial_x u_{low}(t,\cdot) \|_{L^{\infty}(\R)} & \leq  &C |\omega|\, \lambda^{-2-\delta},
\label{vtoro}
\\
\| u_{low}(t,\cdot)- u_{low}(0,\cdot) \|_{L^2(\R)} & \leq  & C |\omega|\, \lambda^{-2-\delta}\, ,
\label{treto}
\end{eqnarray}
if $|t|\leq 1$, $\lambda\gg 1$ and $|\omega|\ll\lambda^{\frac{1-\delta}{2}}$.
\end{lemme} 
\begin{proof}[ Proof of Lemma \ref{ulow}.]
Rescale by setting 
\begin{eqnarray}\label{v}
v(t,x) :=  \lambda^{1+\delta} u_{low}(\lambda^{2+2\delta } t, \lambda^{1+\delta} x).
\end{eqnarray}
Then $v$ is again a solution of the Benjamin-Ono equation.
Since $v(0,x)= -\omega \lambda^{\delta}\tilde \varphi(x)$ we obtain 
$$
\|v(0,\cdot)\|_{H^s}
=
|\omega|\,\lambda^{\delta}\|\tilde{\varphi}\|_{H^s}
$$
and therefore by Proposition~\ref{wp}
\begin{eqnarray}\label{scaling-bis}
\|v(t,\cdot)\|_{H^s}\leq C |\omega|\,\lambda^{\delta},
\end{eqnarray}
if $|t|\leq \min(1, c\,|\omega|^{-4}\lambda^{-4\delta})$
and $s>3/2$. 
But since the right hand-side of (\ref{scaling-bis}) does not
depend on $s$, we conclude that (\ref{scaling-bis}) is valid for any real $s$.
The Sobolev embedding and (\ref{scaling-bis}) now give
\begin{eqnarray}\label{scaling}
\|v_{x}(t,\cdot)\|_{L^{\infty}}\leq C |\omega|\,\lambda^{\delta},
\end{eqnarray}
if $|t|\leq \min(1,c\,|\omega|^{-4}\lambda^{-4\delta})$.

Using (\ref{v}) and the restriction on $|\omega|$, we deduce from (\ref{scaling}) by scaling back that
\begin{eqnarray}\label{lll}
\|\partial_{x}u_{low}(t,\cdot)\|_{L^{\infty}}
\leq C\, |\omega|\,\lambda^{-2-\delta}\, ,
\end{eqnarray}
if $|t|\leq 1$ which is (\ref{vtoro}).

We now turn to the proof of (\ref{parvo}) and (\ref{treto}).
Differentiating (\ref{v}) and using (\ref{scaling-bis}) (with $s=k$) yields
\begin{eqnarray}\label{llll}
\|\partial_{x}^{k}u_{low}(t,\cdot)\|_{L^2}
\leq C |\omega|\,\lambda^{-\frac{1-\delta}{2}-k(1+\delta)},\quad k=0,1,2,\dots
\end{eqnarray}
if $|t|\leq 1$.
Estimate (\ref{llll}) is indeed (\ref{parvo}). Next,
using (\ref{lll}), (\ref{llll}) and the equation satisfied  by $u_{low}$ gives 
$$
\|\partial_{t}u_{low}(t,\cdot)\|_{L^2}
\leq C
\big(
\|\partial_{x}^{2}u_{low}(t,\cdot)\|_{L^2}
+
\|\partial_{x}u_{low}(t,\cdot)\|_{L^{\infty}}\|u_{low}(t,\cdot)\|_{L^2}
\big)
\leq C\, |\omega|\,\lambda^{-2-\delta},
$$
if $|t|\leq 1$.
We now observe that (\ref{treto}) can be deduced from the above bound
via the fundamental theorem of calculus, applied to $u_{low}$ in the time variable.
This completes the proof of Lemma \ref{ulow}. 
\end{proof}
We now set for $\lambda\geq 1$, $0<\delta<1$ and $|\omega| \ll 
\lambda^{\frac{1-\delta}2}$,
\begin{equation}
u_{ap}(t,x):=
u_{low}(t,x)-
\lambda^{-\frac{1}{2}-\frac{\delta}{2}-s}\varphi_\lambda(x)
\cos(-\lambda^2t+\lambda x-\lambda\, t\, u_{low}(0,x)).
\label{ap2}
\end{equation}
The above function is an approximate solution of (\ref{2.2}) for $\lambda\gg 1$ and $s>0$ as shows the next statement.
\begin{lemme}\label{l2}
Let $s>0$, $0<\delta<1$,  $|\omega| \ll
\lambda^{\frac{1-\delta}2}$ and $|t|\le 1$. Set
$$
F:=(\partial_{t}+H\partial_{x}^{2})u_{ap}+u_{ap}\,\partial_{x}u_{ap}.
$$
Then there exist positive constants 
$C$ and $\lambda_0$ such that for 
$\lambda\geq \lambda_0$
one has
$$
\|F(t,\cdot)\|_{L^{2}(\R)}
\leq
C\left( 
\lambda^{-\delta-s} + \lambda^{\frac{1-\delta}2-2s}
\right).
$$ 
\end{lemme}
\begin{proof}
Set $\Phi:=-\lambda^2t+\lambda x+\omega\, t $.
We observe that
\begin{eqnarray*}
u_{ap}(t,x)=u_{low}(t,x)-
\lambda^{-\frac{1}{2}-\frac{\delta}{2}-s}\varphi_{\lambda}(x)\cos \Phi.
\end{eqnarray*}
Furthermore, we define the high frequency part of $u_{ap}$ by setting
$$
u_{h}(t,x):=-\lambda^{-\frac{1}{2}-\frac{\delta}{2}-s}\varphi_{\lambda}(x)\cos \Phi.
$$
Next, we can write 
$$
(\partial_{t}+H\partial_{x}^{2})u_{ap}+u_{ap}\,\partial_{x}u_{ap}
= F_1 + F_2 + F_3 + F_4 +F_5,
$$
where
\begin{eqnarray*}
F_1 & := & (\partial_t + H \partial_{x}^2) u_{low} +  u_{low}\partial_x u_{low}
\\
F_2 & := & 
-\lambda^{-\frac{1+\delta}{2}-s} \,
\cos\Phi\,\, \partial_x \Big( u_{low} \,\varphi_\lambda \Big)
\\
F_3 & := & u_{h} \partial_x u_{h} 
\\
F_4 & := & - \lambda^{-\frac{1+\delta}{2}-s} \Big[H \partial_x^2 ,\varphi_\lambda  \Big] \cos \Phi
\\
F_5 & := & - \lambda^{-\frac{1+\delta}{2}-s}\varphi_\lambda 
(\partial_t + H \partial_{x}^2+u_{low}\partial_x)\cos \Phi \, .
\end{eqnarray*}
Since $u_{low}$ is a solution of the Benjamin-Ono equation, we deduce that $F_1=0$. 
Using that
$\varphi_{\lambda}\tilde{\varphi}_\lambda=\varphi_{\lambda}$, we readily obtain that
$$
F_5  = 
\lambda^{\frac{1-\delta}{2}-s} (u_{low}(t,x) - u_{low}(0,x) )
\varphi_\lambda(x)
\sin \Phi .
$$
Using Lemma \ref{ulow}, we get
\begin{eqnarray}\label{F_5}
\|F_5(t,\cdot)\|_{L^2} \leq C \,
\lambda^{\frac{1-\delta}{2}-s}\,
|\omega| \lambda^{-2-\delta}
\leq C\,\lambda^{-1-2\delta-s}.
\end{eqnarray}
It remains to bound $F_2$, $F_3$ and $F_4$. Let us expand $F_4$ as 
\begin{multline}\label{f4}
F_4 = \lambda^{\frac{3-\delta}2-s}  [H, \varphi_{\lambda}] \cos\Phi
\,+\,2\lambda^{ - \frac{1+\delta}2 - s -\delta}  
H \left((\varphi^{\prime})_\lambda   \sin  \Phi\right) 
- 
\\
-
\lambda^{-\frac52 -\frac{5\delta}2 - s } 
H \left( (\varphi^{\prime\prime})_\lambda  \cos\Phi \right)\, , 
\end{multline}
where
$$
(\varphi^{\prime})_\lambda=\varphi^{\prime}\big(\frac{x}{\lambda^{1+\delta}}\big),\quad
(\varphi^{\prime\prime})_\lambda=\varphi^{\prime\prime}\big(\frac{x}{\lambda^{1+\delta}}\big)\, .
$$
The first term in the right hand-side of (\ref{f4}) is controlled in $L^2$ by the estimate
$$
\|[H, \varphi_{\lambda}] \cos\Phi\|_{L^2}\leq C_{N}\lambda^{-N}
$$
which follows easily from the definition of the Hilbert transform.
The $L^2$ norm of the other terms in the right hand-side of (\ref{f4}) are readily estimated by $ c \lambda^{ - \delta - s}$. 
Therefore
\begin{eqnarray}\label{F4}
\|F_4(t,\cdot)\|_{L^2} \leq C\, \lambda^{-\delta-s}\, .
\end{eqnarray}
Expanding $\partial_{x}u_{h}$, the $L^2$ norm of $F_3$ is controlled as follows 
\begin{eqnarray}\label{F3}
\|F_3(t,\cdot)\|_{L^2}  
\leq
C\,\lambda^{-\frac{3}{2}-\frac{3\delta}{2}-2s}
+C\, \lambda^{\frac{1}{2}-\frac{\delta}{2}-2s}
\leq C\, \lambda^{\frac{1-\delta}{2}-2s} \, .
\end{eqnarray}
Next, using Lemma \ref{ulow} and the assumption on $|\omega|$, we obtain
\begin{multline}\label{F2}
\|F_2(t,\cdot)\|_{L^2}\leq
C\big( 
\lambda^{-s}\|\partial_x u_{low}(t,\cdot)\|_{L^{\infty}}
+
\lambda^{-\frac{3+3\delta}{2}-s}
\|u_{low}(t,\cdot)\|_{L^{2}}
\big)
\leq
\\
\leq C\,\lambda^{-\frac{3+3\delta}{2}-s}.
\end{multline}
Collecting (\ref{F_5}), (\ref{F4}), (\ref{F3}) and (\ref{F2}) completes the proof of Lemma~\ref{l2}.
\end{proof}
Let us now finish the proof of  Theorem~\ref{thm5}.
The first step is to bound $u_{\omega,\lambda}$ in high Sobolev norms. 
We distinguish two cases : $s>3/2$ and $0<s\leq 3/2$.
In the second case we will need to exploit
the higher conservation laws for the Benjamin-Ono equation while in the first case we use Proposition~\ref{wp} instead. 

Let $s>3/2$. Observe that for $3/2<\sigma<s$
$$
\|u_{\omega,\lambda}(0,\cdot)\|_{H^{\sigma}}
\leq
C\big(\lambda^{\sigma-s}+|\omega|\lambda^{-\frac{1-\delta}{2}}\big).
$$
Therefore for $k\geq s$, it follows from Proposition~\ref{wp} that for
$\lambda\gg 1$,
\begin{equation}\label{high1} 
\Vert u_{\omega,\lambda}(t,\cdot) \Vert_{H^k} \leq C \Vert u_{\omega,\lambda}(0,\cdot) \Vert_{H^k}\leq C \lambda^{k-s}, \quad |t|\leq 1. 
\end{equation}
Let $0<s\le \frac32$.
Using the conservation laws associated to the Benjamin-Ono equation
(cf. \cite[Lemma 3.3.2]{ABFS}), we  get the following bound uniformly 
in $t\in\R$ 
\begin{eqnarray}\label{CL}
\|u_{\omega,\lambda}(t,\cdot)\|_{H^{2}}
\leq
C\big(\|u_{\omega,\lambda}(0,\cdot)\|_{H^{2}}+\|u_{\omega,\lambda}(0,\cdot)\|_{L^{2}}^{5}\big)
\leq
C\big(1+  \lambda^{2-s}\big),
\end{eqnarray}
and therefore we obtain 
\begin{equation} \label{high2}
\|u_{\omega,\lambda}(t,\cdot)\|_{H^{2}}
\leq C\lambda^{2-s},\quad t\in\R,\quad \lambda\geq 1\, .
\end{equation}
Let $u_{ap}$ be as in \eqref{ap2}. Set 
$$
v_{\omega,\lambda}:=u_{\omega,\lambda}-u_{ap}\, .
$$
The aim is to show that $v_{\omega,\lambda}$ is small comparing to $u_{ap}$ in the $H^s$ norm.
\\

Due to Lemma \ref{ulow}, we get
$$
\Vert u_{low}(t,\cdot)\Vert_{H^s} 
\leq C\, |\omega|\, \lambda^{-\frac{1-\delta}2},
$$
if $|t|\leq 1$.
Next, using Lemma \ref{loc}, we obtain the bound
$$
\|u_{ap}(t,\cdot)\|_{H^{k}(\R)}
\leq C\, \lambda^{k-s},
$$
if $|t|\leq 1$ and $k\geq s$.

Therefore using (\ref{high1}) and (\ref{high2}), we get the bounds for the  
high Sobolev norms 
\begin{eqnarray} 
\label{firstbound}
\|v_{\omega,\lambda}(t,\cdot)\|_{H^{k}}\leq C\,\lambda^{k-s}, 
\end{eqnarray}
if $|t|\leq 1$ and $3/2<s<k$, and
\begin{eqnarray} 
\label{firstbound-bis}
\|v_{\omega,\lambda}(t,\cdot)\|_{H^{2}}\leq C\lambda^{2-s}, 
\end{eqnarray}
if $t\in\R$ and $0<s<3/2$.
\\

Further, we prove a good bound of the $L^2$ norm of 
$v_{\omega,\lambda}$. Clearly
\begin{equation}\label{Eq-v}
(\partial_{t}+H\partial^{2}_{x})v_{\omega,\lambda}
+
v_{\omega,\lambda}\,\partial_{x}v_{\omega,\lambda}
+
\partial_{x}(u_{ap}\,v_{\omega,\lambda})+
F=0,\quad v_{\omega,\lambda}(0,x)=0
\end{equation}
with 
$$ F = (\partial_t + H \partial_x^2 ) u_{ap} +u_{ap}\partial_x u_{ap}, $$
which satisfies 
$$\|F(t,\cdot)\|_{L^{2}} \leq C\, \lambda^{- \frac{\min \{ \delta, 1-\delta\}}2 -s}$$
by Lemma \ref{l2} and the assumption $1-s<\delta<1$.

The second endpoint in the bounds for $v_{\omega,\lambda}$ is the $L^2$ estimate 
\begin{equation}\label{second}  
\Vert v_{\omega,\lambda}(t,\cdot) \Vert_{L^2} \leq C \lambda^{ - \frac{\min \{ \delta, 1-\delta\}}2 -s },\quad |t|\leq 1. 
\end{equation}  
To prove \eqref{second}, we multiply (\ref{Eq-v}) by $v_{\omega,\lambda}$ and we integrate in $x$ 
\begin{eqnarray*}
\frac{d}{dt}\|v_{\omega,\lambda}(t,\cdot)\|_{L^2}^{2}
\leq 
C\big(\|\partial_{x}u_{ap}(t,\cdot)\|_{L^{\infty}}\|v_{\omega,\lambda}(t,\cdot)\|_{L^2}^{2}+
\|v_{\omega,\lambda}(t,\cdot)\|_{L^2}\|F (t,\cdot)\|_{L^{2}}\big)\, .
\end{eqnarray*}
Hence, since we have for $1-s <\delta<1$ and $\lambda \gg 1$, 
\begin{multline*}
\|\partial_x u_{ap}(t,\cdot)\|_{L^{\infty}}\leq 
C\|\partial_x u_{low}(t,\cdot)\|_{L^{\infty}}+C\lambda^{\frac{1-\delta}2-s} \leq
C|\omega|\lambda^{-2-\delta}+C\lambda^{\frac{1-\delta}2-s} \ll 1,
\end{multline*}
we readily get the bound \eqref{second}.

We now complete the proof by an interpolation argument.
Let first $s>3/2$. 
Choose $k\in [s+\frac{1}{2},s+2]$ and interpolate between
\eqref{firstbound} and \eqref{second} as follows
$$
\Vert v_{\omega,\lambda}(t,\cdot) \Vert_{H^s} 
\le 
\Vert v_{\omega,\lambda}(t,\cdot) \Vert_{L^2}^{\frac{k-s}k} 
\Vert v_{\omega,\lambda}(t,\cdot) \Vert_{H^k}^{\frac{s}k}
\leq C \lambda^{ -\frac{\min\{\delta,1-\delta\}}{4(s+2)}}.   
$$
If $s\le \frac32$ we obtain the same estimate by using $k=2$ in the
interpolation and (\ref{firstbound-bis}) instead of \eqref{firstbound}.
This completes the proof of Theorem~\ref{thm5}.
\end{proof}
We end this section by a series of remarks.
\\

The method of proof of Theorem~\ref{thm5} can be generalized to many other equations.
For example the corresponding to Theorem~\ref{thm5} 
result in the context of the KdV equation provides 
a family of essentially linear KdV waves 
($\omega\rightarrow 0$) as approximate solutions and thus no instability 
property of the flow is displayed.
\\

The proof of Theorem \ref{thm2} is based on a gauge transform reducing (\ref{2.2}) to a problem which,
despite the lack of semi-linear well-posedness displayed by Theorem~\ref{thm3}, 
shares many features with a semi-linearly well-posed problem.
The idea of Tao was further pushed in a series of subsequent papers
\cite{Luc,BP,IK} which enables one to lower the restriction on $s$ in
Theorem~\ref{thm2} and to treat the periodic case.
\\

One may consider the higher dispersion versions of the (\ref{2.2})
\begin{equation}\label{higher}
u_{t}-Lu_{x}+uu_x=0,\quad u(0)=u_0\, ,
\end{equation}
where $L$ is Fourier multiplier with symbol $|\xi|^{\gamma}$, $1\leq \gamma\leq 2$.
The KdV equation corresponds to $\gamma=2$, and, thanks to Theorem~\ref{thm1} in this case (\ref{higher}) is semi-linearly well-posed
in $H^s(\R)$, $s>-3/4$. On the other hand, in view of the result of \cite{MST1}, it seems reasonable to conjecture that for
$1< \gamma< 2$, the Cauchy problem (\ref{higher}) is not semi-linearly well-posed in {\it all} $H^s(\R)$.
\\

Another instance when the notion of semi-linear well-posedness is naturally involved is the analysis of the Cauchy problem
for the Kadomtsev-Petviashvili equations. The  Kadomtsev-Petviashvili (KP) equations are natural two dimensional generalizations
of the KdV equation (cf. \cite{KaPe}). There are two KP equations, the KP-I equation
\begin{equation}\label{KPI}
(u_{t}+u_{xxx}+uu_{x})_{x}-u_{yy}=0,
\end{equation} 
and the KP-II equation
\begin{equation}\label{KPII}
(u_{t}+u_{xxx}+uu_{x})_{x}+u_{yy}=0.
\end{equation} 
It is known that the Cauchy problem for the KP-II equation (\ref{KPII}) is semi-linearly well-posed in $H^{s}(\R^2)$, $s\geq 0$
and even in Sobolev type spaces of negative indices (spaces of distributions), cf. \cite{Bo4,TaTz,IM}.
On the other hand, in view of the result of \cite{MST2}, it seems reasonable to conjecture that the Cauchy problem
for the KP-I equation (\ref{KPI}) is not  semi-linearly well-posed in all $H^s(\R^2)$.
\section{Nonlinear Schr\"odinger equations (NLS)}
\subsection{Nonlinear Schr\"odinger equations on $\R^d$}
Consider the Cauchy problem for the nonlinear Schr\"odinger equation
\begin{equation}\label{3.1}
iu_{t}+\Delta u+|u|^{2}u=0,\quad u(0)=u_{0},
\end{equation}
posed on the Euclidean space $\R^d$, $d\geq 1$. Equation (\ref{3.1}) is a focusing model. The defocusing model 
\begin{equation}\label{3.1-biss}
iu_{t}+\Delta u-|u|^{2}u=0
\end{equation}
is also of interest. The long time dynamics of (\ref{3.1}) and (\ref{3.1-biss}) are quite different.
But for our discussion here (small time analysis) it will be relevant to concentrate only on (\ref{3.1}).
There has been a large number of articles studying (\ref{3.1}), (\ref{3.1-biss}) and their generalizations, when $|u|^{2}u$
(which is the term involved in many applications)
is replaced by a more general nonlinear term $f(|u|^{2})u$ (cf. \cite{GV1,GV2,GV3,Kato2,Kato3,CW,Ts,Y} ...).
The equation (\ref{3.1}) is an infinite dimensional Hamiltonian equation with
canonical coordinates $(u,\bar{u})$ and Hamiltonian
$$
H(u,\bar{u})=\frac{1}{2}\int_{\R^d}|\nabla u|^{2}-\frac{1}{4}\int_{\R^d}|u|^{4}\, .
$$
The Hamiltonian is formally preserved by the flow of (\ref{3.1}). So is the $L^2$ norm of $u$.
Therefore the space $H^1(\R^d)$ is a natural phase space\footnote{In this space the well-posedness of (\ref{3.1-biss}) is actually
global in time. For (\ref{3.1}) the well-posedness is global as far as the data is small in a suitable sense,
for large data solutions developing singularities in finite time appear.}
for (\ref{3.1}) and (\ref{3.1-biss}) at least for $d\leq 4$ when the second term of the
Hamiltonian is dominated by the first one and the $L^2$ norm of $u$.
Fortunately, we can achieve this regularity for $d\leq 3$ in the context of the well-posedness theory of (\ref{3.1}).
More precisely, we have the following result regarding the well-posedness of (\ref{3.1}).
\begin{theoreme}
[cf. \cite{CW}]
\label{thm3.1}
Let $s>\frac{d-2}{2}$, $d\geq 2$. Then the Cauchy problem (\ref{3.1}) is semi-linearly well-posed for data in $H^{s}(\R^{d})$.
\end{theoreme}
\begin{proof}
We will give the proof because it is ``typical'' for a semi-linearly well-posed problem. 
It is worth noticing that such a proof is indeed quite different
from the reasoning in the proof of Proposition~\ref{wp} above. 
To simplify a little the notations we will only consider the case
$d=2$, the proof in higher dimensions being very similar. The proof is based on the following Strichartz 
inequality for the free evolution. 
\begin{proposition}\label{stri}
Let $(p,q)$ such that,
$\frac{1}{p}+\frac{1}{q}=\frac{1}{2}$, $p> 2$.
Then there exists a constant $C>0$ such that for every $T>0$, every $u_{0}\in L^{2}(\R^2)$,
\begin{eqnarray*}
\|e^{it\Delta}u_{0}\|_{L^{p}([0,T]\,;\,L^{q}({\R}^{2}))}\leq C\|u_{0}\|_{L^{2}({\R}^{2})}\, .
\end{eqnarray*}
\end{proposition}
\begin{proof}
The proof of Proposition~\ref{stri} can be found in \cite{Ca}. 
\end{proof}
Let us now show how Proposition~\ref{stri} implies Theorem~\ref{3.1}.
Consider the integral equation corresponding to (\ref{3.1})
\begin{eqnarray}\label{cub}
u(t)=e^{it\Delta}u_0+i\int_{0}^{t}e^{i(t-\tau)\Delta}(|u(\tau)|^{2}u(\tau))d\tau.
\end{eqnarray}
Let us fix a real number $\sigma$ satisfying 
$$
0<\sigma<\min\Big\{s,\frac{1}{2}\Big\} \, .
$$
The value of $\sigma$ being fixed, we define $q\in[2,4[$ by the identity
$$
\frac{1}{q}-\frac{1}{4}=\frac{\sigma}{2}.
$$
Next, we define $p$ such that
$$
\frac{1}{p}+\frac{1}{q}=\frac{1}{2}.
$$
Set
$$
X_{T}=L^{\infty}([0,T]\,;\, H^s(\R^2))\cap Y_T \cap Z_{T},
$$ 
where $Y_T$, $Z_T$ are equipped with the norms
$$
\|u\|_{Y_T}^{2}=\sum_{N-{\rm dyadic}}N^{2s}\|\Delta_{N}(u)\|_{L^{p}_{T}L^{q}(\R^2)}^{2},
\quad
\|u\|_{Z_T}^{2}=\sum_{N-{\rm dyadic}}N^{2s}\|\Delta_{N}(u)\|_{L^{4}_{T}L^{4}(\R^2)}^{2}\,.
$$
The sums over $N$ are running over all dyadic values of $N$, i.e  $N=2^{n}$, $n\geq 0$ and
$$
u=\sum_{N}\Delta_{N}(u)
$$
is a Littlewood-Paley decomposition\footnote{Actually we do not need to use the precise from of $\Delta_N$ here,
cut-off projectors would also make work the argument.} 
of $u$. More precisely, $\Delta_N$ are the Fourier multipliers defined by 
$$
\widehat{\Delta_{N}(u)}(\xi)=\varphi(N^{-1}\xi)\hat{u}(\xi),\quad N=2^{n},\quad n\geq 1
$$
and
$$
\widehat{\Delta_{1}(u)}(\xi)=\psi(\xi)\hat{u}(\xi)
$$
with
$\varphi\in C_{0}^{\infty}(]1,2[)$ and 
$$
\psi(\xi)+\sum_{n\geq 1}\varphi(2^{-n}\xi)=1\, .
$$
Notice that if $u\in {\mathcal S}'(\R^2)$ then $\Delta_{N}(u)$ is localized at frequencies of order $N$.
It is also useful to see the norm in $Y_T$ as $\|N^{s}\Delta_{N}(u)\|_{l^{2}_{N}L^{p}_{T}L^q}$.
\\

Next, using Proposition~\ref{stri}, we get
\begin{eqnarray*}
\|e^{it\Delta}\Delta_{N}(u_0)\|_{L^{p}_{T}L^{q}(\R^2)}+\|e^{it\Delta}\Delta_{N}(u_0)\|_{L^{4}_{T}L^{4}(\R^2)}\leq
C\|\Delta_{N}(u_0)\|_{L^{2}(\R^2)}
\end{eqnarray*}
and therefore
\begin{eqnarray}\label{initiale}
\|e^{it\Delta}u_0\|_{X_T}\leq C\|u_0\|_{H^s}\, .
\end{eqnarray}
Similarly, using the Minkowski inequality, we get the bound,
\begin{equation}\label{23}
\Big\|\int_{0}^{t}e^{i(t-\tau)\Delta}(|u(\tau)|^{2}u(\tau))d\tau
\Big\|_{X_T}\leq C\int_{0}^{T}\Big\|e^{-i\tau\Delta}(|u(\tau)|^{2}u(\tau))\Big\|_{H^s}d\tau.
\end{equation}
In order to bound the right hand-side of (\ref{23}), we will show that there exists $C>0$ such that for every
$u_1$, $u_2$, $u_3$ in $X_T$, 
\begin{eqnarray}\label{sigma}
\int_{0}^{T}\Big\|e^{-i\tau\Delta}(u_{1}(\tau)\overline{u_2}(\tau)u_{3}(\tau))\Big\|_{H^{s}}
d\tau\leq CT^{\sigma}\prod_{j=1}^{3}\|u_j\|_{X_T}.
\end{eqnarray}
By duality, to show (\ref{sigma}), it suffices to obtain that for every $w\in H^{-s}(\R^2)$,
\begin{eqnarray}\label{duality}
\int_{0}^{T}\int_{\R^2}
(u_{1}(\tau)\overline{u_2}(\tau)u_{3}(\tau))(e^{i\tau\Delta}\overline{w})\, d\tau
\leq
CT^{\sigma}\|w\|_{H^{-s}}\prod_{j=1}^{3}\|u_j\|_{X_T}\, .
\end{eqnarray}
Notice that if $u_1(\tau)$, $\overline{u_2}(\tau)$, $u_3(\tau)$ 
are localized at frequencies $N_1$, $N_2$, $N_3$ respectively then only frequencies of order $\leq C(N_1+N_2+N_3)$
of $e^{i\tau\Delta}\overline{w}$ contribute to the left hand-side of (\ref{duality}).
Writhing down the Littlewood-Paley decompositions of
$u_1(\tau)$, $\overline{u_2}(\tau)$, $u_3(\tau)$ et $w$, 
using the  H\"older inequality and Proposition~\ref{stri} to bound  $e^{i\tau\Delta}\overline{w}$,
we deduce that we can bound the left hand-side of (\ref{duality}) by 
\begin{eqnarray}\label{somme}
\sum_{N\leq C(N_1+N_2+N_3)} 
\|\Delta_{N_1} (u_1)\|_{L^{4}_{T}L^{4}}
\|\Delta_{N_2} (u_2)\|_{L^{4}_{T}L^{4}}
\|\Delta_{N_3} (u_3)\|_{L^{4}_{T}L^{4}}
\|\Delta_{N}(w)\|_{L^{2}}
\end{eqnarray}
By symmetry, we can suppose that in (\ref{somme}) the summation is restricted to 
$$
N_1\geq N_2\geq N_3\, .
$$
Next, using the Sobolev embedding  $W^{\sigma,q}(\R^2)\subset L^{4}(\R^2)$, and, the  H\"older inequality in the time 
variable, we get the bound 
\begin{equation*}
\|\Delta_{N_2} (u_2)\|_{L^{4}_{T}L^{4}}\|\Delta_{N_3} (u_2)\|_{L^{4}_{T}L^{4}}
\leq C
T^{\sigma}(N_2 N_3)^{\sigma}\|\Delta_{N_2} (u_2)\|_{L^{p}_{T}L^{q}}\|\Delta_{N_3} (u_3)\|_{L^{p}_{T}L^{q}}\, .
\end{equation*}
Set
$$
c_1(N)=N^{s}\|\Delta_{N_1}( u_1)\|_{L^{4}_{T}L^{4}}\,, \quad
c_{j}(N)=N^{s}\|\Delta_{N} (u_j)\|_{L^{p}_{T}L^{q}}\, , \,j=2,3,\quad
$$
and $d(N)=N^{-s}\|\Delta_{N} (w)\|_{L^{2}}$. 
We obtain that (\ref{somme}) is bounded by
\begin{eqnarray}\label{red}
\sum_{N\leq C N_1}
\sum_{N_2,N_3}
\Big(\frac{N}{N_1}\Big)^{s}\frac{CT^{\sigma}}{(N_2 N_3)^{s-\sigma}}
c_{1}(N_1)c_{2}(N_2)c_{3}(N_3)d(N)\, .
\end{eqnarray}
Summing geometric series in $N_2$, $N_3$, we obtain that (\ref{red}) is bounded by
\begin{eqnarray}\label{29}
CT^{\sigma}
\|u_2\|_{X_T}\|u_3\|_{X_T}\,
\sum_{N\leq  C N_1}
\Big(\frac{N}{N_1}\Big)^{s}
c_{1}(N_1)d(N)\, .
\end{eqnarray}
To bound (\ref{29}), we use the following lemma.
\begin{lemme}\label{trik}
For every $\Lambda>0$, every $s>0$ there exists $C>0$ such that
if $(c_{N_{0}})$ and $(d_{N_1})$ are two sequences of nonnegative numbers
indexed by the dyadic integers, then,
\begin{equation*}
\sum_{N_{0}\leq \Lambda N_{1}}\, 
\Big(\frac{N}{N_1}\Big)^{s}\,c_{N_{0}}\,d_{N_1}
\leq C\Big(\sum_{N_0}c_{N_{0}}^{2}\Big)^{\frac{1}{2}}\Big(\sum_{N_1}d_{N_1}^{2}\Big)^{\frac{1}{2}}\, .
\end{equation*}
\end{lemme}
\begin{proof}
Let us set
$$
K(N_0,N_1):=
\11_{N_{0}\leq \Lambda N_{1}}\,\frac{N_{0}^{s}}{N_1^{s}}\, .
$$
Summing geometric series imply that there exists $C>0$ such that
$$
\sup_{N_0}\sum_{N_1}K(N_0,N_1)+\sup_{N_1}\sum_{N_0}K(N_0,N_1)\leq C\, .
$$
Therefore the Schur lemma implies the boundedness on $l^{2}_{N_0}\times l^{2}_{N_1}$  of the bilinear form with 
kernel $K(N_0,N_1)$. This completes the proof of Lemma \ref{trik}.
\end{proof}
Using Lemma~\ref{trik}, we bound (\ref{29}) 
by the right hand-side if (\ref{duality}) which completes the proof of (\ref{sigma}).

Estimates  (\ref{initiale}), (\ref{23}) et (\ref{sigma}) yield that for every bounded set $B$ of $H^{s}(\R^2)$
there exists $T>0$ such that for every  $u_{0}\in B$ the right hand-side of (\ref{cub}) is a contraction in a suitable 
ball of $X_T$.

Let us finally explain how we obtain the propagation of regularity property for data in $H^{\tilde{s}}$, $\tilde{s}\geq s$.
Denote by $X_{T}^{s}$ the space $X_T$ used above, associated to the $H^s$ regularity.
It is easy to observe that the preceding analysis also gives the tame estimate
\begin{eqnarray*}
\Big\|
\int_{0}^{t}e^{i(t-\tau)\Delta}(|u(\tau)|^{2}u(\tau))d\tau
\Big\|_{X_T^{\tilde{s}}}\leq CT^{\sigma}\|u\|_{X_{T}^{s}}^{2}\|u\|_{X_{T}^{\tilde{s}}}
\end{eqnarray*}
which implies the propagation of the $H^{\tilde{s}}$ regularity in a straightforward way.
This ends the discussion on the proof of Theorem~\ref{thm3.1}.
\end{proof}
The indice $\frac{d-2}{2}$ appeared in Theorem~\ref{thm3.1} is closely related to the scaling of the equation (\ref{3.1}). More precisely
if $u(t,x)$ solves (\ref{3.1}) then so does
$$
u_{\lambda}(t,x):=\lambda u(\lambda^{2}t,\lambda x).
$$
The norm of $u_{\lambda}$ in the homogeneous Sobolev $\dot{H}^{s}$ is independent of $\lambda$ only for $s=\frac{d-2}{2}$.

At this point, it is worth noticing that the scaling invariance is responsible for the existence of solutions of (\ref{3.1})
which concentrate in a point. Such kind of concentrations may give ill-posedness results only below the scaling norm.
As we will see later concentration on higher dimensional objects as curves are responsible for ill-posedness above
the scaling exponent. 
\\

It turns out that the result of Theorem~\ref{thm3.1} is essentially sharp, i.e. the point concentrations coming from of the scaling
invariance are the worst possible.
\begin{theoreme}
[cf. \cite{CCT2}]
\label{thm3.2}
Let $d\geq 2$. Then :
\begin{enumerate}
\item
For $d=2$, the Cauchy problem (\ref{3.1}) is not semi-linearly well-posed for data in $H^{s}(\R^d)$,
$s<0\, (=\frac{d-2}{2})$. Moreover, it is not well-posed for data in $H^{s}(\R^d)$, $s\leq -1\, (=-\frac{d}{2})$.
\item
For $d\geq 3$, the Cauchy problem (\ref{3.1}) is not well-posed for data in $H^{s}(\R^d)$,
$0<s<\frac{d-2}{2}$ or $s< -\frac{d}{2}$. Moreover, it is not semi-linearly well-posed for data in $H^{s}(\R^d)$,
$-\frac{d}{2}\leq s\leq 0$.
\end{enumerate}
\end{theoreme}
\begin{proof}
In order to simplify the exposition, we will give the proof of  Theorem~\ref{thm3.2}, for $d\geq 5$ and $s$ a positive
integer. This will cover the most interesting case $s=1$, i.e. the ill-posedness of (\ref{3.1}) in the ``energy space'' 
$H^1(\R^d)$, $d\geq 5$. The proof of  Theorem~\ref{thm3.2} in the other cases has a very similar flavor.
\\

Let us first observe that it suffices to prove the following statement.
\begin{proposition}\label{CCT}
Let $d\geq 5$ and $s\in ]0,\frac{d-2}{2}[$ be a positive integer. 
Then there exist a sequence $(t_n)$ of positive numbers tending to zero and a sequence 
$(u_{n}(t))$ of $C^{\infty}(\R^d)\cap H^{s}(\R^d)$ functions defined for $t\in [0,t_n]$, such that 
$$
(i\partial_t+\Delta) u_n +|u_n|^2 u_n=0
$$
with
$$
\lim_{n\rightarrow \infty}\|u_{n}(0,\cdot)\|_{H^s(\R^d)}=0\, ,
$$
and
$$
\lim_{n\rightarrow \infty}\|u_{n}(t_n,\cdot)\|_{H^s(\R^d)}= \infty \,.
$$
\end{proposition}
\begin{proof}
Let us consider an initial data concentrating in the point $x=0$
$$
u_{n}(0,x) := \kappa_n n^{\frac{d}{2}-s} \varphi(n x),\quad n\gg 1,
$$ 
where $\varphi$ is a non identically zero smooth compactly supported function and 
$$
\kappa_n= \log^{- \delta_1}(n)
$$ 
with $\delta_1>0$ to be fixed later. Remark that 
$$
\|u_{n}(0,\cdot)\|_{H^s(\R^d)} \sim \kappa_n.
$$
The function 
$$
v_{n}(t,x)=\kappa_n n^{\frac{d}{2}-s}\varphi(n x) e^{it [\kappa_n n^{\frac{d}{2}-s}\varphi(n x)]^{2}}
$$ 
is the solution of the equation
\begin{equation}\label{novo}
i\partial_t  v_n +|v_n|^{2}v_n=0,\quad v_n(0,x)= u_{n}(0,x)\, .
\end{equation}
It turns out that for very small times $v_n$ is near the actual 
solutions\footnote{A similar idea was used, in a different context, by
Kuksin \cite{Kuksin}.} of (\ref{3.1}). 
\\

Next, for a fixed integer $l>d/2$, we define quantity,
$$
E_n(u):= \Big(n^{2s} \|u\|^2_{L^2(\R^d)} + n^{-2(l-s)}\|u\|^2_{H^{l}(\R^d)}\Big)^{\frac{1}{2}}
$$
which can be seen as a semi-classical energy of $u$. Notice that, uniformly in $n$, 
\begin{equation}\label{control}
\|u\|_{H^s}\leq CE_{n}(u)\, .
\end{equation}
The main point in the proof of Proposition~\ref{CCT} is the next statement.
\begin{lemme}\label{le.est}
Fix $\delta_2\in\R$ such that
$$
0<\delta_2<\frac{1}{l+1}\, .
$$
Then the solution $u_n$ of \eqref{3.1} with initial data 
$$
u_{0}(x)=\kappa_n n^{\frac{d}{2}-s} \varphi(n x)
$$ 
exists for $0\leq t\leq t_n$, where 
$$
t_n = \log^{\delta_2}(n)n^{-2(\frac{d}{2}-s)}.
$$ 
Moreover, there exists $\varepsilon >0$ such that for $t\in[0,t_n]$,
$$
E_n(u_n(t)- v_n(t))\leq C n^{-\varepsilon}\, .
$$
\end{lemme}
\begin{proof}
Since the initial data are in $H^{l}$, $l>d/2$, we know that $u_{n}(t)$ exist on small time interval $[0,\widetilde{t}_n]$.
Consequently, to prove Lemma~\ref{le.est}, we simply prove the {\em a priori} estimates which ensure, by a classical
bootstrap argument, both the existence and the control on $E_n(u_n(t)- v_n(t))$ for $t\in[0,t_n]$.
Let us set 
$$
w_n:=u_n -v_n\, .
$$
The a  priori estimates involved in the proof are energy inequalities in the equation satisfied by $w_n$,
\begin{eqnarray*}
(i \partial_t + \Delta)w_n
&  = & -\Delta v_n-v_{n}^{2}\overline{w}_n-2|v_n|^{2}w_n-2v_n |w_n|^{2}-\overline{v}_n w_n^2-|w_n|^{2}w_n
\\
& := & -\Delta v_n +\Lambda(v_n,w_n)\, .
\end{eqnarray*}
Using the explicit formula for $v_n$, we have that for $0\leq t \leq t_n$,
\begin{equation}\label{vn1}
\| v_n(t,\cdot)\|_{L^\infty(\R^d)} \leq C n^{\frac{d}{2}-s} 
\end{equation}
and for $\sigma\geq 0$, $0\leq t \leq t_n$,
\begin{equation}\label{vn2}
\| v_n(t,\cdot)\|_{H^{\sigma}(\R^d)} \leq C n^{\sigma-s}  \log^{\delta_2\sigma}(n)\, .
\end{equation}
Let us now estimate $E_{n}(\Lambda(v_{n}(t),w_{n}(t)))$ for $0\leq t \leq t_n$. 
Using the Gagliardo-Nirenberg inequality
$$
\|f\|_{L^{\infty}(\R^d)}\leq C\|f\|_{L^{2}(\R^d)}^{1-\frac{d}{2l}}\|f\|_{H^{l}(\R^d)}^{\frac{d}{2l}},\quad l>d/2,
$$
we infer that
\begin{equation}\label{GN}
\|f\|_{L^{\infty}(\R^d)}\leq Cn^{\frac{d}{2}-s}E_{n}(f)\, .
\end{equation}
Coming back to the expression for $\Lambda(v_n,w_n)$, we get
$$
\|\Lambda(v_n,w_n)\|_{L^2}\leq C\big(\|v_n\|_{L^{\infty}}^{2}+\|w_n\|_{L^{\infty}}^{2}\big)\|w_n\|_{L^{2}}\, .
$$
Using (\ref{vn1}) and (\ref{GN}), we obtain that for $0\leq t \leq t_n$,
$$
n^{s}\|\Lambda(v_n(t),w_n(t))\|_{L^2}\leq Cn^{2(\frac{d}{2}-s)}\big(  E_n(w_n(t))+E_{n}^{3}(w_n(t))\big)\, .
$$
Next, using several times the classical bilinear inequality
$$
\|fg\|_{H^l}\leq C\big(\|f\|_{L^{\infty}}\|g\|_{H^{l}}+ \|g\|_{L^{\infty}}\|f\|_{H^{l}} \big),
$$
we deduce that
\begin{multline*}
\|\Lambda(v_n,w_n)\|_{H^l}\leq C\big(
\|v_n\|_{L^{\infty}}^{2}\|w_n\|_{H^l}
+
\|v_n\|_{L^{\infty}}\|v_n\|_{H^l}\|w_n\|_{L^{\infty}}
+
\\
+
\|v_n\|_{H^l}\|w_n\|_{L^{\infty}}^{2}
+
\|v_n\|_{L^{\infty}}\|w_n\|_{L^{\infty}}\|w_n\|_{H^l}
+
\|w_n\|_{L^{\infty}}^{2}\|w_n\|_{H^l}
\big)\, .
\end{multline*}
Using (\ref{vn1}), (\ref{vn2}) and (\ref{GN}), we infer that  for $0\leq t \leq t_n$,
$$
n^{-(l-s)}\|\Lambda(v_n(t),w_n(t))\|_{H^l}
\leq Cn^{2(\frac{d}{2}-s)}
\log^{\delta_2 l}(n)
\big(  E_n(w_n(t))+E_{n}^{3}(w_n(t))\big)\, .
$$
Summarizing the above discussion yields that for $0\leq t\leq t_n$,
$$
E_{n}(\Lambda(v_n(t),w_n(t)))
\leq Cn^{2(\frac{d}{2}-s)}
\log^{\delta_2 l}(n)
\big(  E_n(w_n(t))+E_{n}^{3}(w_n(t))\big)\, .
$$
Next, we estimate the source term $-\Delta v_n$. Using (\ref{vn2}), we obtain that for $0\leq t \leq t_n$,
$$
\|\Delta v_{n}(t,\cdot)\|_{L^2}\leq C\|v_{n}(t,\cdot)\|_{H^2}\leq Cn^{2-s}\log^{2\delta_2 }(n)\, .
$$
Similarly,
$$
\|\Delta v_{n}(t,\cdot)\|_{H^l}\leq C\|v_{n}(t,\cdot)\|_{H^{l+2}}\leq Cn^{l+2-s}\log^{(l+2)\delta_2 }(n)\, .
$$
Therefore, for $0\leq t\leq t_n$,
$$
E_{n}(\Delta v_n(t))\leq Cn^{2}\log^{(l+2)\delta_2 }(n)\, .
$$
Coming back to the equation solved by $w_n$, we deduce that for $0\leq t\leq t_n$,
\begin{multline*}
\frac d {dt} E^{2}_{n}(w_n(t)) 
\leq  C\, n^{2(\frac{d}{2}-s)} \log^{\delta_2 l}(n)\big(E_n^{2} (w_n(t))+E_n^{4}(w_n(t))\big)+
\\
+
Cn^2 \log^{(l+2)\delta_2}(n)E_n(w_n(t))  \, .
\end{multline*}
Suppose first that $E_{n}(w_n(t)) \leq 1$ which is clearly the case at least for $t\ll 1$ since $w_n(0,x)=0$.
Using the elementary inequality
$$
2n^2 \log^{(l+2)\delta_2}(n)E_n(w_n)
\leq  n^{2(\frac{d}{2}-s)} \log^{\delta_2 l}(n)E_n^{2}(w_n)
+
n^{4-2(\frac{d}{2}-s)}\log^{(l+4)\delta_2}(n)
$$
we obtain that for $0\leq t\leq t_n$,
\begin{equation*}
\frac d {dt} 
\Big[e^{-C\,t\,  n^{2(\frac{d}{2}-s)} \log^{\delta_2 l}(n)}E^{2}_{n}(w_n(t)) \Big]
\leq  C\, n^{4- 2(\frac{d}{2}-s) }\log^{(l+4)\delta_2}(n)\,e^{-C\,t\,  n^{2(\frac{d}{2}-s)} \log^{\delta_2 l}(n)} \, .
\end{equation*}
Integrating between $0$ and $t$ yields,
$$
E_{n}(w_n(t))\leq Cn^{2-2(\frac{d}{2}-s)}\, \log^{2\delta_2}(n)\,e^{C\log^{(l+1)\delta_2}(n)}\, .
$$
The assumption $s<(d-2)/2$ implies $2-2(\frac{d}{2}-s)<0$. Moreover $\delta_2$ is such that $(l+1)\delta_2<1$.
Therefore there exists $\varepsilon>0$ such that
$$
E_{n}(w_n(t))\leq C n^{-\varepsilon}\, .
$$
Finally, a bootstrap argument allows to drop the assumption 
$$
E_{n}(w_n(t)) \leq 1\, .
$$
This completes the proof of Lemma~\ref{le.est}.
\end{proof}
Let us now finish the proof of  Proposition~\ref{CCT}.
We need to make a proper choice of the number $\delta_1$ involved in the definition of $\kappa_n$.
Using the explicit formula for $v_n$, we easily obtain that for $s\in\Z^{+}$,
$$
\|v_{n}(t_n,\cdot)\|_{H^s}\geq C\kappa_{n}\big(t_n[\kappa_n n^{\frac{d}{2}-s}]^{2}\big)^{s},
$$
provided $t_n[\kappa_n n^{\frac{d}{2}-s}]^{2}\gg 1$, i.e. $\log^{\delta_2-2\delta_1}(n)\gg 1$.
The first assumption on $\delta_1$ is thus
$$
0<\delta_1<\frac{\delta_2}{2}\,.
$$
Therefore, using Lemma~\ref{le.est} and (\ref{control}), we obtain that for $n\gg 1$,
\begin{multline*}
\|u_{n}(t_n,\cdot)\|_{H^s}\geq 
C\|v_{n}(t_n,\cdot)\|_{H^s}-Cn^{-\varepsilon}
\geq 
C\kappa_{n}\big(t_n[\kappa_n n^{\frac{d}{2}-s}]^{2}\big)^{s}-Cn^{-\varepsilon}
=
\\
=
C\log^{s\delta_2-(1+2s)\delta_1}(n)-Cn^{-\varepsilon}\, .
\end{multline*}
The proof of Proposition~\ref{CCT} is completed by choosing $\delta_1\in\R$ such that
$$
0<\delta_1<\frac{s\delta_{2}}{1+2s}\, .
$$
\end{proof}
The ansatz with $v_n$ as an approximate solution still holds for $-\frac{d}{2}<s\leq 0$ and very small times of order
$\sim \log^{\delta}(n)n^{-2(\frac{d}{2}-s)}$ with a suitable $\delta>0$. 
Unfortunately, we can no longer bound from below $\|v_{n}(t_n,\cdot)\|_{H^s}$ as above. For that reason 
one can not get the failure of well-posedness for data in $H^{s}(\R^d)$, $-\frac{d}{2}<s\leq 0$.
One can however still obtain the lack of semi-linear well-posedness, by an argument very similar to the one that we
presented in the context of the Benjamin-Ono equation.
This ends the discussion on the proof of Theorem~\ref{thm3.2}.
\end{proof}
\begin{remarque}\label{rrr}
The approach to Theorem~\ref{thm3.2} we present here may seem more involved than in \cite{CCT2} but it has the advantage to avoid
the scaling considerations of \cite{CCT2}. In particular it works for variable coefficients second order operators instead of $\Delta$
or for (\ref{3.1}) posed on a curved space.
\end{remarque}
Let us next consider (\ref{2.1}) posed on the real line $\R$. In this case the critical threshold for the semi-linear well-posedness
is shifted with respect to the scaling regularity.
\begin{theoreme}
[cf. \cite{CW,CCT1,Ts}]
\label{thm3.3}
For $d=1$ the Cauchy problem (\ref{3.1}) is semi-linearly well-posed for data in $H^s(\R)$, $s\geq 0$, and, 
it is not semi-linearly well-posed for data in $H^s(\R)$, $s<0$. 
\end{theoreme}
It is worth noticing that in the proof of the lack of semi-linear well-posedness for $s<0$, one uses a family of solutions
which concentrate on the line $\{t=0\}$ of the space time $(t,x)$. 
This family of solutions is related to the Galilean invariance of~(\ref{3.1}).
\\

Notice that for $d\leq 4$, the space $H^1$ is essentially covered by
Theorem~\ref{thm3.1}. For $d\geq 5$, the result of Theorem~\ref{thm3.1} is far from the regularity $H^1$,
and moreover as we have shown, for $d\geq 5$ the Cauchy problem (\ref{3.1}) is not well-posed for data in $H^1(\R^d)$.
\\

In order to have an $H^1$ theory in dimensions $d\geq 5$, it is reasonable to replace (\ref{3.1}) with the equation
\begin{eqnarray}\label{gen}
(i\partial_{t}+\Delta)u=F(u),
\end{eqnarray}
where the nonlinear interaction $F$ is supposed to satisfy $F(0)=0$ and is supposed of the form 
$
F=\frac{\partial V}{\partial \bar{z}}
$
with a {\it positive} $V\in C^{\infty}(\C\,;\,\R)$ satisfying 
$V(e^{i\theta}z)=V(z)$, $\theta\in\R$, $z\in\C$, and, for some $\alpha>1$,
\begin{equation*}
|\partial_{z}^{k_1}\partial_{\bar{z}}^{k_2}\, V(z)|\leq C_{k_1,k_2}(1+|z|)^{1+\alpha-k_1-k_2}\, .
\end{equation*}
The number $\alpha$ involved in the second condition on $V$ corresponds to the degree of the nonlinearity $F(u)$.
The Hamiltonian associated to (\ref{gen}) is
\begin{equation}\label{intr-3}
\int_{\R^d} |\nabla u|^{2}+\int_{\R^d}V(u)\, .
\end{equation}
which controls the $\dot{H}^{1}$ norm.
If $\alpha<1+\frac{4}{d-2}$, the second term in (\ref{intr-3}) is controlled by the first one and the $L^2$ norm of $u$. 
It is therefore reasonable to expect that for $\alpha<1+\frac{2}{d-2}$ the Cauchy problem for (\ref{gen}) is well-posed for data
in $H^1$. It turns out to be the case at least for $d=5,6$.
\begin{theoreme}
[cf. \cite{GV1,Kato1,BGT4}]
\label{thm3.4}
Let $\alpha<1+\frac{4}{d-2}$ and $d\leq 6$. Then  the Cauchy problem associated to  (\ref{gen}) 
is semi-linearly well-posed for data in $H^1(\R^d)$.
\end{theoreme}
\begin{remarque}
The proof of Theorem~\ref{thm3.4} relies crucially on the Strichartz inequalities for the non-homogeneous linear problem
(recall that in the proof of Theorem~\ref{thm3.1} we only used the Strichartz inequalities for the free evolution).
For $d\geq 7$ one can still prove the existence and the uniqueness.
Since the nonlinearity in (\ref{gen}) is not polynomial,
the propagation of of the regularity is a nontrivial problem in the analysis of (\ref{gen})  (cf. e.g. \cite{Gi-DEA}).
In addition, even if one can prove the regularity propagation, the semi-linear well-posedness of (\ref{gen}) for $d\geq 7$
remains a non trivial issue. 
\end{remarque}
\subsection{Nonlinear Schr\"odinger equations on compact manifolds}
Let $(M,g)$ be a compact smooth boundaryless Riemannian manifold of dimension $d\geq 2$. Denote by $\Delta_g$ the 
Laplace operator associated to the metric $g$. In this section we consider the analog of (\ref{3.1-biss}) on $(M,g)$
\begin{equation}\label{3.1-bis}
iu_{t}+\Delta_{g} u-|u|^{2}u=0,\quad u(0)=u_{0}\, .
\end{equation}
We have the following well-posedness result for (\ref{3.1-bis})
\begin{theoreme}
[cf. \cite{BGT1}]
\label{thm3.5}
The Cauchy problem (\ref{3.1-bis}) is semi-linearly well-posed for data in $H^{s}(M)$, $s>\frac{d-1}{2}$.
\end{theoreme}
The proof of Theorem~\ref{thm3.5} is based on the following Strichartz inequality with derivative losses for the free evolution
$$
\|\exp(it\Delta_{g})u_0\|_{L^{p}([0,T]\,;\,L^{q}(M))}\leq C_{T}\|u_{0}\|_{H^{\frac{1}{p}}(M)},
$$
where 
\begin{equation}\label{admissible}
\frac{2}{p}+\frac{d}{q}=\frac{d}{2},\quad p\geq 2,\quad (p,q)\neq (2,\infty)\, .
\end{equation}
Then the contraction argument is performed in the space
$$
X_{T}=L^{\infty}([0,T]\,;\,H^s(M))\cap L^{p}([0,T]\,;\, W^{s-\frac{1}{p},q}(M)),
$$
with suitable $(p,q)$ satisfying (\ref{admissible}).

Notice that there is a gap between the regularity $\frac{d-1}{2}$ of  Theorem~\ref{thm3.5} and the
regularity $\frac{d-2}{2}$ of Theorem~\ref{3.1}. It is a natural question what happens for data in $H^s$,
$s\in [\frac{d-2}{2},\frac{d-1}{2}]$. It turns out that in the case of the flat torus $\T^d$ we can recover the result
of the case $\R^d$. 
\begin{theoreme}
[cf. \cite{Bo1}]
\label{thm3.6}
Let $M$ be the flat torus. Then (\ref{3.1-bis}) is semi-linearly well-posed for data in $H^{s}(\T^d)$, $s>\frac{d-2}{2}$.
\end{theoreme}
As in the case of the KdV equation of the previous section,
the proof of Theorem~\ref{thm3.6} uses the Fourier transform restrictions spaces. 
An important additional element in the analysis is the use of bilinear improvements of the Strichartz inequalities.

The spaces of Bourgain and bilinear Strichartz estimates can also be used in the case of the sphere to get the following result.
\begin{theoreme}
[cf. \cite{BGT5,BGT6}]
\label{thm3.7}
Let $M$ be the standard sphere $S^d$. Then : 
\begin{enumerate}
\item
If $d=2$ then (\ref{3.1-bis}) is semi-linearly well-posed for data in $H^{s}(S^2)$, $s>\frac{1}{4}$.
\item
If $d\geq 3$ then (\ref{3.1-bis}) is semi-linearly well-posed for data in $H^{s}(S^d)$, $s>\frac{d-2}{2}$.
\end{enumerate}
\end{theoreme}
It is worth noticing that the ill-posedness argument of Theorem~\ref{thm3.2} still applies in the setting of the Riemannian manifolds
(cf. Remark \ref{rrr}).
Therefore, for $d\geq 3$, the indice $s=\frac{d-2}{2}$, turns out to be the critical one for both the well-posedness and the
semi-linear well-posedness of (\ref{3.1-bis}) posed on the flat torus or on the standard sphere. 

We observe that in the case $d=2$ the assumption $s>1/4$ in Theorem~\ref{thm3.7} is more restrictive then in the case of the torus $\T^2$.
It turns out that this assumption is sharp, as far as the semi-linear well-posedness is concerned.
\begin{theoreme}
[cf. \cite{BGT2}]
\label{thm3.8}
Let $M$ be the standard sphere $S^2$. Then  (\ref{3.1-bis}) is not semi-linearly well-posed for data in $H^{s}(S^2)$, 
$s<\frac{1}{4}$.
\end{theoreme}
\begin{proof}
The main ingredient in Theorem~\ref{thm3.8} is a description of the evolution by the flow of  (\ref{3.1-bis})
of the highest weight spherical harmonics.
\begin{proposition}\label{highest}
Let $T>0$, $s\in ]\frac{3}{20},\frac{1}{4}[$, $\kappa\in ]0,1[$.
Take $M=S^2$ with the canonical metric in (\ref{3.1-bis}).
For $n\in \N$, we denote by $\psi_{n} : S^2\rightarrow \C$ the restriction to $S^2$ of the harmonic polynomial $(x_{1}+ix_{2})^n$. 
Then the solution $u_{n}(t)$ of (\ref{3.1-bis}) with initial data 
$\kappa\varphi_n$, where $\varphi_n=n^{\frac{1}{4}-s}\psi_n$ is globally defined, and, for $t\in [0,T]$ it can be represented as
\begin{eqnarray}\label{pot}
u_n(t)=\kappa\,
e^{-it(n(n+1)+\kappa^{2}\omega_{n})}\,\big(\varphi_n+r_n(t)\big),
\end{eqnarray}
where $\omega_n\approx n^{\frac{1}{2}-2s}$
and $r_n(t)$ satisfies
\begin{equation*}
\|r_n(t)\|_{H^{s}(S^2)}\leq
C_{T}\, n^{-\delta}
\end{equation*}
where $\delta>0$ and $C_{T}$ depends on $T$ but not on $n$.
Moreover there exists $C>0$, independent of $T$ and $n$ such that
\begin{equation}
\|u_n\|_{L^{\infty}(\R\,;\,H^{s}(S^2))}\leq C\kappa.
\end{equation}
\end{proposition}
\begin{proof}[Proof of Proposition~\ref{highest}]
Recall that $\psi_n$ is and eigenfunction of $-\Delta_g$ associated to and eigenvalue $n^2+n$.
An easy computation shows that 
$$
\|\psi_n\|_{L^{p}}\approx n^{-\frac{1}{2p}},\quad n\gg 1\, . 
$$
Therefore $\|\varphi_n\|_{H^s}\approx 1$ and $\|\varphi_{n}\|_{L^2}\leq Cn^{-s}$.
Similarly to the Euclidean case, the solutions of  (\ref{3.1-bis}) enjoy the conservation laws
\begin{equation}\label{cons1}
\int_{M}|u(t,x)|^{2}dx={\rm Const}
\end{equation}
and
\begin{equation}\label{cons2}
\int_{M}|\nabla u(t,x)|^{2}dx+\frac{1}{2}\int_{M}|u(t,x)|^{4}dx={\rm Const}\, .
\end{equation}
The $H^1$ well-posedness result of Theorem~\ref{thm3.5} applies for the initial data $u_0=\kappa\varphi_n$ and
we obtain a local solution $u_n(t)$.  
Using the conservation laws (\ref{cons1}), (\ref{cons2}), we deduce that  
the $H^1$ norm of $u_{n}(t)$ is uniformly bounded.
Therefore, we can reiterate the well-posedness results and to obtain that the solutions $u_{n}(t)$ are globally defined.
\\

For every $\alpha\in\R$, we denote by $R_{\alpha}$ the rotation of $\R^{3}$ defined by
$$
R_{\alpha}(x_1,x_2,x_3)=(\cos\alpha\, x_1-\sin\alpha\, x_2,\, \sin\alpha\, x_1+\cos\alpha\, x_2,\, x_3)
$$
and by $R_{\alpha}^{\star}$ the associated unitary operator of $L^{2}(S^2)$,
$$
R_{\alpha}^{\star}u(x)=u(R_{\alpha}(x)).
$$
Observe that $R_{\alpha}^{\star}\psi_{n}=\exp(in\alpha)\psi_{n}$ for every $\alpha\in\R$.
The following elementary lemma will be useful in the sequel.
\begin{lemme}\label{A}
Let $n\in \Z^{+}$ and $u\in L^{2}(S^2)$ be such that for every $\alpha\in\R$,
\begin{eqnarray}\label{pg1}
R_{\alpha}^{\star}u=\exp(in\alpha)u.
\end{eqnarray}
Then the decomposition of $u$ in spherical harmonics reads
$$
u=c\psi_n+\sum_{j}g_{j}
$$
where $c\in\C$ and each $g_{j}$ is a spherical harmonic of degree $>n$.
\end{lemme}
\begin{proof}
Since the family $(R_{\alpha}^{\star})_{\alpha\in\R}$ is a one-parameter group of unitary operators leaving invariant 
the space of spherical harmonics of degree $l$, one can find an orthonormal basis $(h_k)$ of 
$L^{2}(S^2)$ such that, for every $k$, $h_k$ is a spherical harmonic satisfying, for some $n_k\in\Z$, for every $\alpha\in\R$,
\begin{eqnarray}\label{pg2}
R_{\alpha}^{\star}h_k =\exp(in_{k}\alpha)h_{k}.
\end{eqnarray}
Comparing (\ref{pg1}) and (\ref{pg2}), the decomposition of 
$u$ in the basis $(h_k)$ reads \begin{eqnarray}\label{pg3}
u=\sum_{k\,:\,n_{k}=n}c_{k}h_{k}.
\end{eqnarray}
Let $h$ be a spherical harmonic of degree $l$ satisfying property (\ref{pg1}) for every $\alpha\in\R$.
Denote by $P$ the $l$-homogeneous polynomial on $\R^{3}$
such that $h=P_{|S^2}$. Then (\ref{pg1}) is equivalent to
\begin{eqnarray}\label{pg4}
\forall\, x\in\R^{3},\quad P(R_{\alpha}(x))=\exp(in\alpha)P(x).
\end{eqnarray}
Let us decompose $P$ according to the powers of $z=x_{1}+ix_{2}$
and $\bar{z}=x_{1}-ix_{2}$
\begin{eqnarray}\label{pg5}
P(x_1,x_2,x_3)
=
\sum_{p+q\leq l}a_{pq}\, z^{p}\bar{z}^{q}x_3^{l-p-q}
\end{eqnarray}
where $a_{pq}\in\C$. In view of (\ref{pg4}), (\ref{pg5}) and
\begin{eqnarray}\label{pg6}
P(R_{\alpha}(x))
=
\sum_{p+q\leq l}a_{pq}\exp(i(p-q)\alpha)
z^{p}\bar{z}^{q}x_3^{l-p-q},
\end{eqnarray}
we conclude that $a_{pq}=0$ unless $p-q=n$. As a consequence,
$$
l\geq p+q\geq p-q=n
$$
and, if $l=n$, then $p=n$ and $q=0$, so that $P=cz^{n}$, i.e. $h=c\psi_n$ for some $c\in\C$.
Coming back to decomposition (\ref{pg3}) completes the proof of Lemma~\ref{A}.
\end{proof}
Using Lemma~\ref{A}, we can write
$$
|\varphi_n|^{2}\varphi_n=\omega_n \varphi_n+r_{n}\, ,
$$
where $r_{n}$ contains only spherical harmonics of degree $>n$ in its spectral decomposition and
$$
\omega_{n}=\frac{\|\varphi_n\|_{L^4}^{4}}{\|\varphi_n\|_{L^2}^{2}}\approx n^{\frac{1}{2}-2s}\, .
$$
Observe that that $R_{\alpha}^{\star}u_{n}$ is a solution of (\ref{3.1-bis}) with data $u_0=\kappa e^{in\alpha}\varphi_n$.
On the other hand $e^{in\alpha}u_{n}$ is also a solution of (\ref{3.1-bis}) with the same initial data.
Therefore, using the uniqueness assertion of Theorem~\ref{thm5} (in spaces invariant under the action of $R_{\alpha}^{\star}$) 
for the Cauchy problem  (\ref{3.1-bis}), we obtain 
$$
R_{\alpha}^{\star}u_{n}=e^{in\alpha}u_{n}\, .
$$
Using Lemma~\ref{A}, we deduce that $u_{n}(t)$ is a linear combination of $\psi_n$ and spherical harmonics of degree $>n$.
\\

Let us give the heuristic argument which permits us to find an ansatz for $u_{n}(t)$. In view of the above discussion, 
we may hope that $u_{n}(t)$ can be written as
$$
u_{n}(t)=\kappa c_{n}(t)\varphi_n+ {\rm ``small\,\, error'' }\, .
$$
Substituting this in the equation (\ref{3.1-bis}), neglecting the ``small error'' and projection on $\varphi_n$
yields the equation
$$
i\kappa \dot{c_n}-n(n+1)\kappa c_n-\kappa^{3}\,\omega_{n}\,|c_n|^{2}c_{n}=0,\quad c_{n}(0)=1
$$
which gives
$$
c_{n}(t)=e^{-it(n(n+1)+\kappa^{2}\omega_n)}\, .
$$

In order to make the above formal discussion rigorous, we set
$$
u_{n}(t)=\kappa\, e^{-it(n(n+1)+\kappa^{2}\omega_n)}((1+z_{n}(t))\varphi_n+q_{n}(t)),
$$
where $z_{n}(0)=0$, $q_{n}(0)=0$ and $q_{n}(t)$ contains only spherical harmonics of degree $>n$ in its spectral decomposition.
Proposition~\ref{highest} is clearly a consequence of the following statements.
\begin{lemme}\label{lem3}
There exists a constant $C>0$, independent of $T$ and $n$ such that
$$
\|q_n(t)\|_{H^s}\leq Cn^{-\frac{1}{4}-s}.
$$
\end{lemme}
\begin{lemme}\label{lem4}
There exists a constant $C_{T}>0$, which depends on $T$ 
but not on $n$ such that
$$
\sup_{t\in[0,T]}
|z_n(t)|\leq C_{T}n^{\frac{1}{4}-3s}.
$$
\end{lemme}
\begin{proof}[Proof of Lemma~\ref{lem3}]
Let us first rewrite the conservation laws (\ref{cons1}), (\ref{cons2}) in terms of $z_n(t)$ and $q_n(t)$. 
Since $\varphi_n$ is orthogonal to $q_n(t)$ in $L^{2}(S^2)$ as well as $\nabla \varphi_n$ to $\nabla q_n(t)$,
we can rewrite (\ref{cons1}) and (\ref{cons2}) as
\begin{equation}\label{cons1'}
|1+z_n(t)|^{2}\|\varphi_n\|_{L^2}^{2}+\|q_n(t)\|_{L^2}^{2}=\|\varphi_n\|_{L^2}^{2},
\end{equation}
\begin{multline}\label{cons2'}
|1+z_n(t)|^{2}\|\nabla\varphi_n\|_{L^2}^{2}+\|\nabla q_n(t)\|_{L^2}^{2}+\frac{1}{2\kappa^2}\|u_n(t)\|_{L^4}^{4}
=
\\
\|\nabla\varphi_n\|_{L^2}^{2}+\frac{\kappa^2}{2}\|\varphi_n\|_{L^4}^{4}.
\end{multline}
Observe that
$$
\|\nabla\varphi_n\|_{L^2}^{2}=\langle\varphi_n\,,\,-\Delta\,\varphi_n\rangle=n(n+1)\|\varphi_n\|_{L^2}^{2},
$$
where $\langle\cdot,\cdot\rangle$ denotes the $L^{2}(S^2)$ scalar product.
Therefore multiplying (\ref{cons1'}) with $-n(n+1)$ and adding it to
(\ref{cons2'}) gives
\begin{equation}\label{star}
\|\nabla q_n(t)\|_{L^2}^{2}-n(n+1)\|q_n(t)\|_{L^2}^{2} \leq
\frac{\kappa^2}{2}\|\varphi_n\|_{L^4}^{4}
\leq C n^{\frac{1}{2}-4s}\,  .
\end{equation}
We can decompose
$$
q_n(t)=\sum_{l\geq n+1}q_{n,l}(t)
$$ 
where $q_{n,l}\in\mbox{Ker}(\Delta_{S^2}+l(l+1))$.
Hence
\begin{equation*}
\|\nabla q_n(t)\|_{L^2}^{2}-n(n+1)\|q_n(t)\|_{L^2}^{2}
=
\sum_{l\geq n+1}
(l(l+1)-n(n+1))\|q_{n,l}(t)\|_{L^2}^{2}.
\end{equation*}
If $l\geq n+1$, we have $l(l+1)-n(n+1)\geq l$ and therefore
$$
\|\nabla q_n(t)\|_{L^2}^{2}-n(n+1)\|q_{n}(t)\|_{L^2}^{2}\geq 
\|q_{n}(t)\|_{H^{\frac{1}{2}}}^{2}\, .
$$
Coming back to (\ref{star}),
\begin{equation}\label{edno}
\|q_{n}(t)\|_{H^{\frac{1}{2}}}\leq Cn^{\frac{1}{4}-2s}\, .
\end{equation}
On the other hand, we also have that if $l\geq n+1$ then $l(l+1)-n(n+1)\geq n$ and thus
$$
n\|q_n(t)\|_{L^2}^{2}\leq 
Cn^{\frac{1}{2}-4s}
$$
which implies
\begin{equation}\label{dve}
\|q_{n}(t)\|_{L^2}\leq Cn^{-\frac{1}{4}-2s}\, .
\end{equation}
Using (\ref{edno}) and (\ref{dve}), we finally arrive at
$$
\|q_n(t)\|_{H^s}\leq \|q_{n}(t)\|_{L^2}^{1-2s}\|q_{n}(t)\|_{H^{\frac{1}{2}}}^{2s}
\leq Cn^{-\frac{1}{4}-s}\, .
$$
This completes the proof of Lemma \ref{lem3}
\end{proof}
\begin{proof}[Proof of Lemma~\ref{lem4}]
Let us set $w_{n}(t):=z_{n}(t)\varphi_n+q_{n}(t)$. By projecting the equation
$$
(i\partial_t+\Delta)u_n-|u_n|^{2}u_n=0
$$
on the mode $\varphi_n$, we get that $z_n$ solves the equation
\begin{multline*}
i\dot{z}_{n}+\kappa^{2}\omega_n z_n= 
\frac{\kappa^2}{\|\varphi_n\|_{L^2}^2}
\Big(\langle|\varphi_n+w_n|^{2}(\varphi_n+w_n)\, , \,\varphi_n\rangle-\langle |\varphi_n|^{2}\varphi_n\, , \,\varphi_n\rangle\Big)
=
\\
= 
\frac{\kappa^2}{\|\varphi_n\|_{L^2}^2}
\int\,\Big( (2|\varphi_n|^{2}w_n+\varphi_n^{2}\overline{w_n})\overline{\varphi_n}
+2\Re(\overline{\varphi_n}w_n)w_n\overline{\varphi_n}+|w_n|^{2}|\varphi_n|^{2}
+|w_n|^{2}w_n\overline{\varphi_n}\Big).
\end{multline*}
Next, the equation for $z_n(t)$ can be rewritten as
\begin{multline*}
i\dot{z}_{n}+\kappa^{2}\omega_n z_n = 
2\omega_n\kappa^{2}z_n+\omega_n\kappa^{2}\bar{z_n} 
+
\frac{\kappa^2}{\|\varphi_n\|_{L^2}^2}
{\mathcal O}\Big(|z_n|^{2}\int|\varphi_n|^{4}
+
\\
+
|z_n|^{3}\int|\varphi_n|^{4}
+\int|q_n|^{3}|\varphi_n|
+\int|q_n|^{2}|\varphi_n|^{2}+|\langle q_n,r_n\rangle|\Big).
\end{multline*}
Let us estimate the source terms. Write using Lemma~\ref{lem3}
\begin{multline*}
\frac{\int|q_n|^{3}|\varphi_n|}{\|\varphi_n\|_{L^2}^2}\leq 
Cn^{2s}\|q_n\|_{L^3}^{3}\|\varphi_n\|_{L^{\infty}}
\leq 
Cn^{2s}\|q_n\|_{L^2}\| q_n\|_{H^{\frac{1}{2}}}^{2}\|\varphi_n\|_{L^{\infty}}
\leq 
\\
\leq
Cn^{2s}\,n^{-\frac{1}{4}-2s}\,n^{\frac{1}{2}-4s}\,n^{\frac{1}{4}-s}= 
Cn^{\frac{1}{2}-5s}.
\end{multline*}
Further we have
\begin{equation*}
\frac{\int|q_n|^{2}|\varphi_n|^{2}}{\|\varphi_n\|_{L^2}^2}\leq Cn^{2s}\|q_n\|_{L^2}^{2}\|\varphi_n\|_{L^{\infty}}^{2}
\leq Cn^{2s}\,n^{-\frac{1}{2}-4s}\,n^{\frac{1}{2}-2s}=  Cn^{-4s}.
\end{equation*}
and 
\begin{equation*}
\frac{|\langle q_n\, , \,r_n\rangle|}{\|\varphi_n\|_{L^2}^2}\leq 
Cn^{2s}\|q_n\|_{L^2}\|r_n\|_{L^{2}}
\leq Cn^{2s}\,n^{-\frac{1}{4}-2s}\,n^{\frac{1}{2}-3s}= 
Cn^{\frac{1}{4}-3s}.
\end{equation*}
Therefore, if $s>\frac{1}{8}$, the equation for $z_n(t)$ can be written as
\begin{equation}\label{15}
i\partial_{t} z_{n}=2\omega_n\kappa^{2}\Re(z_n)+{\mathcal O}(\omega_n|z_n|^{2}+\omega_n|z_n|^{3}
+n^{\frac{1}{4}-3s})
\end{equation}
with $z_n(0)=0$. 
Moreover using once again the $L^2$ conservation law (\ref{cons1'}), we have
$$
1-|1+z_n|^{2}=
\frac{\|q_n(t)\|_{L^2}^{2}}{\|\varphi_n\|_{L^2}^{2}}
={\mathcal O}(n^{-\frac{1}{2}-2s}).
$$
Therefore
$$
|2\Re(z_n)+|z_n|^{2}|={\mathcal O}(n^{-\frac{1}{2}-2s})
$$
and the equation (\ref{15}) takes the form
$$
i\partial_{t} z_{n}={\mathcal O}(\omega_n|z_n|^{2}+\omega_n|z_n|^{3}
+n^{\frac{1}{4}-3s}),
$$
with $\omega_n={\mathcal O}(n^{\frac{1}{2}-2s})$.
Hence if we set
$$
M_n(T)=\sup_{0\leq t\leq T}|z_n(t)|,
$$
we obtain
\begin{eqnarray}\label{16}
M_n(T)\leq CT(n^{\frac{1}{2}-2s}[M_n(T)]^{2}+n^{\frac{1}{2}-2s}[M_n(T)]^{3}+n^{\frac{1}{4}-3s}).
\end{eqnarray}
In view of (\ref{16}), we set
$$
\tilde{M_n}(T)=n^{3s-\frac{1}{4}}M_n(T)
$$
and therefore (\ref{16}) yields
$$
\tilde{M_n}(T)\leq CT(1+n^{\frac{3}{4}-5s}[\tilde{M_n}(T)]^{2}+n^{1-8s}[\tilde{M_n}(T)]^{3}).
$$
Since $\tilde{M_n}(0)=0$ and $s>\frac{3}{20}$, we obtain that 
$\tilde{M_n}(T)\leq CT$ uniformly with respect to $n$.
This completes the proof of Lemma \ref{lem4}.
\end{proof}
This proof of  Proposition~\ref{highest} is now completed.
\end{proof}
Notice that the assertion of Theorem~\ref{thm3.8} is particular for the sphere
$S^2$ only for $0\leq s<1/4$. 
Indeed, for $s<0$ we can apply the argument of Theorem~\ref{thm3.2} in the
context of an arbitrary Riemannian manifold.
\\

Let us now show how Proposition~\ref{highest} implies Theorem~\ref{thm3.8} for $3/20<s<1/4$.
The main point is that for $s<1/4$ we have $\omega_n\rightarrow\infty$ as $n\rightarrow\infty$. 
Let us fix $T>0$, $\kappa\in ]0,1]$ and let $(\kappa_n)$ be a sequence of positive numbers such that
$$
(\kappa^{2}-\kappa_{n}^{2})\omega_n=n^{\beta},\quad 0<\beta<1\, .
$$
Since $s<1/4$ and $\omega_n\approx n^{\frac{1}{2}-2s}$, we have that for $\beta\ll 1$, $\kappa_{n}\rightarrow\kappa$.
Let $(u_{\kappa,n})$ and $(u_{\kappa_{n},n})$ be the solutions of (\ref{3.1-bis}) with data $\kappa\varphi_n$ and
$\kappa_n\varphi_{n}$ respectively. Then
$$
\|u_{\kappa,n}(0,\cdot)-u_{\kappa_{n},n}(0,\cdot)\|_{H^s}\leq C|\kappa-\kappa_n|\longrightarrow 0
$$
but thanks to Proposition~\ref{highest}, for $t\in [0,T]$,
$$
\|u_{\kappa,n}(t,\cdot)-u_{\kappa_{n},n}(t,\cdot)\|_{H^s}\geq c|e^{itn^{\beta}}-1|-C_{T}n^{-\delta}
$$
with $\delta>0$. 
The proof of Theorem~\ref{thm3.8} for $3/20<s<1/4$  is completes by observing that for all $n\gg 1$,
$$
\sup_{0\leq t\leq T}|e^{itn^{\beta}}-1|=2\, .
$$

When $1/8< s\leq 3/20$, we need to perform a slight modification of the argument.
Indeed, in this case is suffices to remark that in fact we need to justify the ansatz only on a a {\it small interval}\footnote{Recall 
that a similar idea idea was used in the discussion around Theorem~\ref{thm3.2}.}
$[0,T_n]$ with $T_n$ satisfying
$$
\lim_{n\rightarrow\infty}n^{\frac{1}{2}-2s}T_n=\infty\, .
$$
The bound of Lemma~\ref{lem3} is uniform in time. We only need to slightly modify the proof of Lemma~\ref{lem4}. 
In the case $1/8< s\leq 3/20$, we define 
$$
\tilde{M_n}(T):=T_{n}^{-1}n^{3s-\frac{1}{4}}M_n(T),\qquad 0\leq T\leq T_{n}\, ,
$$
and, the argument of Lemma~\ref{lem4} yields the bound
$$
\tilde{M_n}(T)\leq
C\big(1+T_{n}^{2}n^{\frac{3}{4}-5s}[\tilde{M_n}(T)]^{2}+T_{n}^{3}n^{1-8s}[\tilde{M_n}(T)]^{3}\big)\,.
$$
If we set
$$
T_{n}:=n^{\frac{5s}{2}-\frac{3}{8}-\varepsilon}
$$
with $0< \varepsilon\ll 1$, a bootstrap argument gives
$$
\tilde{M_n}(T)\leq C
$$
and therefore
$$
|z_{n}(t)|\leq Cn^{\frac{1}{4}-3s}\, n^{\frac{5s}{2}-\frac{3}{8}-\varepsilon}=Cn^{-\frac{1}{8}-\frac{s}{2}-\varepsilon},\qquad
0\leq t\leq T_n\, .
$$
Therefore the ansatz is valid on $[0,T_n]$, which is a sufficiently large small
interval to get the instability property of the flow map.
Indeed
$$
n^{\frac{1}{2}-2s}T_n=n^{\frac{1}{8}+\frac{s}{2}-\varepsilon}
$$
which gives the big oscillations needed to assure the lack of uniform continuity of the flow map.
For case $0\leq s\leq 1/8$ we refer to the work of Banica \cite{Banica}, 
where the ansatz of Proposition~~\ref{highest}
is justified up to time one for $0< s\leq 3/20$.
This completes the discussion on the proof of  Theorem~\ref{thm3.8}.
\end{proof}
We end this section by several remarks.
\\

The result of Theorem~\ref{thm3.8} is another instance when we see that the critical indice for the semi-linear well-posedness
is shifted from the scaling one because of concentration on a curve (a closed geodesic). 
It would be interesting to develop a notion of critical
exponent associated to curve similarly to the one associated to a point via the scaling invariance.
\\

It would be interesting to decide whether for some $0\leq s\leq 1/4$, the Cauchy problem (\ref{3.1-bis}) is well-posed
(probably after a suitable gauge transform) for data in $H^s(S^2)$. Recall that such a phenomenon is not excluded as shows
the experience with the modified KdV equation.
\\

We do not know for an analog of Theorem~\ref{thm3.4} in the setting of compact manifolds. 
Moreover, it is known that in case of the sphere $S^6$ the assertion of Theorem~\ref{thm3.4} fails.
More precisely the Cauchy problem
$$
iu_{t}+\Delta_{g} u=(1+|u|^{2})^{\frac{\alpha}{2}}u,\quad u(0)=u_0,\quad 0<\alpha\leq 1,
$$
posed on $S^6$ is not semi-linearly well-posed for data in $H^{1}(S^6)$ (cf. \cite{BGT2,BGT4}).
\\

We refer to the work \cite{BGT3} where the approach of 
Theorem~\ref{thm3.8} is extended to (\ref{3.1-bis}) posed on the unit disc of $\R^2$ with Dirichlet boundary conditions.
We also refer to \cite{CCT3} for ill-posedness results for the cubic NLS posed on the circle $S^1$.
\\

In all our examples for the failure of well-posedness or semi-linear
well-posedness, the leading part of the approximate solutions is on the high
frequencies. We refer to \cite{BT,CCT3} for examples when the main part of the
approximate solution is on the low frequencies (after a high-high interaction).
\section{Final remarks}
There has been a number of works, closely related to the discussion in these notes for nonlinear wave equations
(cf. \cite{BK,CCT2,DG,Li,Leb1,Leb2} ... ).  In the context of the  nonlinear wave equations, again, {\it families} 
of solutions concentrating
at a point contradict the well-posedness (or semi-linear well-posedness) below the scaling exponent. The finite propagation speed
of the wave equation is exploited in \cite{Leb2} to construct a {\it single} solution, concentrating in an infinite number of points,
which stays bounded in $H^s$ (for some suitable $s$) and becomes instantaneously very large in $H^{\sigma}$, $\sigma>s$.
It would be interesting to prove the analogue of Lebeau's result in the context of the NLS. Despite the lack of the finite
propagation speed for the Schr\"odinger operator, the reasoning in the proof of Theorem~\ref{thm3.2} above is of a semi-classical nature
(cf. also \cite{BGT1}) and thus finite propagation speed considerations could be employed. It is worth noticing that, again, 
in the case of nonlinear wave equations ill-posedness above the scaling is closely related to concentrations on curves,
for instance the Lorentz invariance provides families of solutions concentrating on light rays.

The problematic discussed in these notes fits naturally in the context of parabolic PDE's. There has been some first results
in that direction (cf. \cite{Canone, MRY} ...), and, we believe there is further progress to come.

\end{document}